\documentclass[a4paper,11pt]{article}
\usepackage{amsmath,amssymb,amsthm,a4wide,color}
\usepackage[utf8]{inputenc}
\usepackage[T1]{fontenc}
\usepackage{authblk,mathrsfs,empheq}
\usepackage{algorithm,algpseudocode}
\usepackage{graphicx,subcaption}

\usepackage{booktabs}
\usepackage{ragged2e}

\usepackage{caption}

\usepackage{pdfsync}
\newcommand{\RR}{\mathbb{R}}
\newcommand{\CC}{\mathbb{C}}
\newcommand{\mrm}{\mathrm}
\newcommand{\mA}{\mathrm{A}}
\newcommand{\mB}{\mathrm{B}}
\newcommand{\mL}{\mathrm{L}}
\newcommand{\mH}{\mathrm{H}}
\newcommand{\mJ}{\mathrm{J}}
\newcommand{\mV}{\mathrm{V}}

\newcommand{\Id}{\mathrm{Id}}

\newcommand{\mT}{\mathrm{T}}
\newcommand{\mS}{\mathrm{S}}

\newcommand{\mR}{\mathrm{R}}
\newcommand{\mQ}{\mathrm{Q}}
\newcommand{\mX}{\mathrm{X}}
\newcommand{\mbH}{\mathbb{H}}

\newcommand{\mbX}{\mathbb{X}}

\newcommand{\infsup}{\mathop{\mrm{infsup}}}

\newcommand{\bell}{\boldsymbol{\ell}}
\newcommand{\bp}{\boldsymbol{p}}
\newcommand{\bff}{\boldsymbol{f}}

\newcommand{\bx}{\boldsymbol{x}}

\newcommand{\bn}{\boldsymbol{n}}
\newcommand{\bu}{\boldsymbol{u}}
\newcommand{\bv}{\boldsymbol{v}}
\newcommand{\bw}{\boldsymbol{w}}
\newcommand{\bq}{\boldsymbol{q}}

\newcommand{\bg}{\boldsymbol{g}}

\newcommand{\dir}{\textsc{d}}
\newcommand{\neu}{\textsc{n}}

\newcommand{\lbr}{\lbrack}
\newcommand{\rbr}{\rbrack}

\newcommand{\ctru}{\mathfrak{u}}
\newcommand{\ctrv}{\mathfrak{v}}
\newcommand{\ctrw}{\mathfrak{w}}

\newcommand{\supp}{\mathrm{supp}}

\newtheorem{defn}{Definition}[section]
\newtheorem{lem}[defn]{Lemma}
\newtheorem{thm}[defn]{Theorem}
\newtheorem{prop}[defn]{Proposition}
\newtheorem{cor}[defn]{Corollary}
\newtheorem{example}[defn]{Example}

\newcommand*\widefbox[1]{\framebox[.65\columnwidth]{#1}}

\title{\textbf{Non-local optimized Schwarz method }\\
  \textbf{with physical boundaries}}
\author[1]{X.Claeys}

\affil[1]{Sorbonne Université,
  Laboratoire Jacques-Louis Lions}

\date{}

\begin{document}

\maketitle

\begin{abstract}
  We extend the theoretical framework of non-local optimized Schwarz
  methods as introduced in [Claeys,2021], considering an Helmholtz equation
  posed in a bounded cavity supplemented with a variety of conditions
  modeling material boundaries. The problem is reformulated equivalently as
  an equation posed on the skeleton of a non-overlapping partition of the
  computational domain, involving an operator of the form "identity +
  contraction". The analysis covers the possibility of resonance
  phenomena where the Helmholtz problem is not uniquely solvable. In
  case of unique solvability, the skeleton formulation is proved
  coercive, and an explicit bound for the coercivity constant is provided
  in terms of the inf-sup constant of the primary Helmholtz boundary value problem.
\end{abstract}

\section*{Introduction}

Large scale simulation of harmonic wave propagation phenomena remains
a challenge in the context of which one of the most effective
substructuring domain decomposition methods (DDM) was introduced by
Després \cite{MR1291197}. Commonly referred to as
Optimized Schwarz Method (OSM), it consists in local solves of the wave
equation, maintaining a coupling between subdomains
through a reformulation of transmission conditions in terms of ingoing
and outgoing Robin traces. The new transmission conditions involve an
exchange operator that swaps traces from both sides of each interface
between neighboring subdomains. This approach was put in a general
theoretical framework in \cite{MR1764190} and we point to
\cite{zbMATH07020343} for an overview of this type of strategy.

In a discrete setting, the appropriate definition of the exchange
operator raises issues at cross-points, where at least three degrees of
freedom have to communicate, because it is then unclear what should be
the discrete counterpart of swapping. Although several
heuristics had been proposed in the literature for dealing with this
situation \cite{Gander2013,MR3519297,Modave2020b,MR4480644,Bendali2006}, most
strategies based on this local swapping operator experienced
deteriorated performance in the presence of cross points.

\quad\\
In a series of articles \cite{claeys2019new,claeys2021nonself,MR4507159,MR4433119},
we proposed a variant of OSM where the
usual local swapping exchange operator is replaced by an alternative a priori
non-local operator that naturally accommodates the presence of cross-points. This new approach
can cope with arbitrary subdomain partitions, with a possibly very complicated wire basket.
In \cite{claeys2019new}, we analyzed this new approach at the continuous level considering
a transmission problem posed on the full space $\RR^d$, and the formulation associated
to this new DDM strategy was proved strongly coercive, which paved the way to
convergence estimates for linear solvers (e.g. Richardson, GMRes).

This novel approach was adapted to a finite element discretised setting
and a full convergence theory was developed in \cite{MR4433119,claeys2021nonself}. In passing,
this new theoretical framework covered the case of the original Després algorithm
hence offering a genuine generalization. The whole theory was confirmed by
numerical results both in 2D and 3D. While the previous developments were concerned
with scalar harmonic wave propagation, the case of Maxwell's equations was considered
in \cite{MR4507159,parolin:tel-03118712}.

\quad\\
In the present contribution we extend the theory of \cite{claeys2019new} in several
directions.  First of all, while \cite{claeys2019new} considered only the case of
a transmission problem posed on the whole of $\RR^d$, we consider here the case of
a cavity problem posed in a bounded domain $\Omega\subset\RR^d$. This boundary value
problem takes the form 
\begin{equation}\label{PbInit1}
  \begin{aligned}
    & \mrm{div}(\mu^{-1}\nabla u) + \kappa^{2}u = -f\;\text{in}\;\Omega\\
    & + \textrm{boundary condition on $\partial\Omega$.}    
  \end{aligned}
\end{equation}
Here again we reformulate it  as an equation in terms of traces posed on the skeleton
of the subdomain partition, which we call skeleton formulation. While in previous
contributions the problem had been assumed  uniquely solvable (see e.g.
\cite[\S 1]{MR4433119} or \cite[\S 1.2]{claeys2021nonself}),
the analysis is here extended so as to cover the case where \eqref{PbInit1}
is not necessarily uniquely solvable which covers the case of non-trivial resonance
phenomenon. The skeleton formulation is then proved uniquely solvable if and only
if this holds for \eqref{PbInit1} and, if this condition is fulfilled, the skeleton
formulation is proved to be strongly coercive. Although coercivity was already established
in \cite{claeys2019new}, we provide in addition an explicit estimate of the coercivity
constant in terms of the inf-sup condition of the primary variational formulation.

Our whole analysis rests on an interpretation of the properties of \eqref{PbInit1}
in terms of a pair of two closed linear manifolds: one that models transmission conditions,
and another one that models local wave equations. Studying properties of operators by means
of pairs of closed linear manifolds follows the spirit of \cite[\textsc{iv}.4 \& \textsc{iv}.5]{MR1335452}.

Like \cite{claeys2019new}, the present contribution is purely theoretical. It aims at laying solid
analytical foundations for a better understanding of the spectral properties of the skeleton
formulation, which is important in the perspective of devising both computationally efficient
eigensolvers and domain decomposition preconditionners. We do not provide any numerical experiment.
Such results shall be presented in a forthcoming contribution that will develop
a discrete variant of the present analysis, in the spirit of \cite{MR4433119,claeys2021nonself}.

\quad\\
The outline of this article is as follows. In the first two sections we introduce general notations
for both Hilbert analysis and Sobolev spaces, including trace operators, Dirichlet-to-Neumann maps
and harmonic liftings. Next we describe the problem under study, specifying precisely the assumptions
underlying our analysis, which allows in particular to deal with a variety of boundary conditions.
How to apply this framework for common boundary conditions is illustrated with examples. Further notations
are introduced  for dealing with multi-domain configurations. This leads in particular to a
characterization of transmission conditions based on a non-local exchange operator, see
Proposition \ref{NonlocalTrCond}, which had been an important innovation of \cite{claeys2019new}.
We use this multi-domain formalism to re-express the boundary value problem under study. The kernel
and the range of this operator are then re-interpreted in terms of
a pair of closed linear manifolds. One manifold models wave equations local to each subdomain,
and the other one models transmission conditions. Wave equations local to each subdomain are then
re-expressed by means of a so-called scattering operator, which we use to finally provide
a formulation involving tuples of Robin traces on the skeleton of the subdomain partition.
This skeleton formulation is proved to systematically admit closed range, and its kernel is put
in correspondence with the kernel of the original formulation. Finally we prove strong coercivity
for the skeleton formulation and derive an estimate for the coercivity constant that is explicit
with respect to the inf-sup constant of the original variational formulation.

\section{General notation conventions}\label{NotationConventions}

We first set a few general notation conventions regarding
analysis in Banach spaces. All vector spaces that we are going
to consider have $\mathbb{C}$ as scalar field.
Assuming that $\mH$ is a Banach space equipped with the norm
$\Vert \cdot\Vert_{\mH}$, its topological dual denoted $\mH^{*}$ will
systematically be equipped with the norm
\begin{equation}\label{DualNorm}
  \Vert \varphi\Vert_{\mH^{*}} =
  \sup_{v\in \mH\setminus\{0\}} \frac{\vert \varphi(v) \vert}{\Vert v\Vert_{\mH}}.
\end{equation}
The canonical duality pairing will be systematically denoted
$\langle \cdot,\cdot\rangle:\mH^{*}\times\mH\to \mathbb{C}$ and
defined by $\langle \varphi, v\rangle:= \varphi(v)$. Although
the space $\mH$ does not appear explicitly in the notation
"$\langle \varphi, v\rangle$", when such pairing angle brackets
are used, it shall be clear from the context which pair of spaces
$(\mH,\mH^{*})$ is under consideration. We emphasize
that the duality pairings we consider do not involve any
complex conjugation. We shall write $\langle v,\varphi\rangle =
\langle \varphi, v\rangle\;\forall v\in\mH,\varphi\in \mH^*$ indifferently.
For any subset $\mX\subset \mH$, we denote its polar set by
\begin{equation}\label{DefPolarSet}
  \mX^{\circ} :=\{ \varphi\in \mH^*,\langle \varphi,v\rangle = 0\;\forall v\in \mX\}. 
\end{equation}
Assuming that $\mV$ is another Banach space
equipped with the norm $\Vert \cdot\Vert_{\mV}$,
and $\mL:\mH\to \mV$ is a bounded linear map, we shall
refer to its inf-sup constant denoted and defined as follows
\begin{equation}\label{DefInfSupCst}
  \infsup_{\mH\to \mV}(\mL):=\inf_{u\in \mH\setminus\{0\}}
  \frac{\Vert \mL(u)\Vert_{\mV}}{\Vert u\Vert_{\mH}}
\end{equation}
In the case where $\mL$ is invertible, this inf-sup constant
equals the inverse to the continuity modulus of $\mL^{-1}$.
The inf-sup constant is well defined even if $\mL$ is not
invertible though. The adjoint to the map $\mL:\mH\to \mV$ shall be
defined as the unique bounded linear map $\mL^*:\mV^{*}\to \mH^{*}$
satisfying
\begin{equation}\label{DefAdjoint}
  \langle \mL^{*}(p),u\rangle := 
  \langle p,\mL(u)\rangle 
\end{equation}
for all $p\in \mV^*$ and all $u\in\mH$. Once again, we insist that
no complex conjugation comes into play in \eqref{DefAdjoint}. The bounded
linear map $\mL$ induces another bounded linear map
$\overline{\mL}:\mH\to \mV$ defined by $\overline{\mL}(\overline{u}):=\overline{\mL(u)}$
for all $u\in\mH$.

\quad\\
A bounded linear operator $\mT:\mH\to \mH^*$
is called self-adjoint if $\overline{\mT} = \mT^{*}$ and, in this case
we have $\langle \mT(u),\overline{u}\rangle\in \RR$ for all $u\in\mH$.
It is called positive definite if $\langle \mT(u),\overline{u}\rangle
\in (0,+\infty)$ for all $u\in\mH\setminus\{0\}$. If $\mT$ is
both self-adjoint and positive definite, the sesquilinear form
$u,v\mapsto \langle \mT(u),\overline{v}\rangle$ induces a scalar
product over $\mH$ and the associated norm is denoted
\begin{equation}\label{OperatorInducedNorm}
  \Vert u\Vert_{\mT}:=\sqrt{\langle \mT(u),\overline{u}\rangle}.
\end{equation}
We shall also consider cartesian products $\mH_{1}\times \dots\times \mH_{\mJ}$
where each $\mH_{j}$ is a Banach space equipped with the norm $\Vert\cdot\Vert_{\mH_{j}}$.
Then the cartesian product shall be equipped with the following canonical
norm and duality pairings
\begin{equation}\label{CartesianProductNorm}
  \begin{aligned}
    & \Vert \bv\Vert_{\mH_{1}\times \dots\times \mH_{\mJ}}^{2} :=
    \Vert v_1\Vert_{\mH_{1}}^{2}+\dots +\Vert v_{\mJ}\Vert_{\mH_{\mJ}}^{2}\\
    & \langle \bv,\bq \rangle :=
    \langle v_1,q_1\rangle +\dots  +    \langle v_\mJ,q_\mJ\rangle.
  \end{aligned}
\end{equation}
for $\bv = (v_1,\dots,v_\mJ), v_j\in\mH_{j}$, and
$\bq = (q_1,\dots,q_{\mJ}), q_{j}\in\mH^{*}_{j}$.
If $\mV_{j},j=1,\dots,\mJ$ is another collection of Banach spaces
and $\mL_{j}:\mH_{j}\to \mV_{j}$ are bounded linear maps, we
shall also consider the block-diagonal operator
$\mrm{diag}(\mL_{1},\dots,\mL_{\mJ})$,
mapping $\mH_{1}\times \dots\times \mH_{\mJ}$ into
$\mV_{1}\times \dots\times \mV_{\mJ}$ and defined,
for $\bv = (v_1,\dots,v_\mJ), v_j\in\mH_{j}$, and
$\bq = (q_1,\dots,q_{\mJ}), q_j\in\mV_{j}^*$, by
\begin{equation*}
  \langle \bq, \mrm{diag}(\mL_{1},\dots,\mL_{\mJ})\,\bv\rangle :=
  \langle q_1,\mL_{1}(v_1)\rangle + \dots +
  \langle q_{\mJ},\mL_{\mJ}(v_{\mJ}) \rangle.
\end{equation*}

\section{Single domain functional setting}\label{SingleDomainVF}

Now we need to introduce classical function spaces.
For any Lipschitz open set $\omega\subset \RR^d$, we consider
$\mL^{2}(\omega):=\{ v:\omega\to \CC \;\text{measurable},\;
\Vert v\Vert_{\mL^{2}(\omega)}^{2}:=\int_{\omega}
\vert v(\bx)\vert^{2}d\bx<+\infty\}$ and define
Sobolev spaces 
\begin{equation}\label{H1fctspace}
  \begin{aligned}
    & \mH^{1}(\omega):=\{v\in\mL^{2}(\omega),\;
    \nabla v\in\mL^{2}(\omega)^{d} \}\\
    & \Vert v\Vert_{\mH^{1}(\omega)}^{2}:=
    \Vert \nabla v\Vert_{\mL^{2}(\omega)}^{2} + \gamma^{-2}\Vert v\Vert_{\mL^{2}(\omega)}^{2}
  \end{aligned}
  \end{equation}
where $\gamma>0$ is a real positive parameter.
Incorporating $\gamma$-dependency in the norm will allow
to establish $\gamma$-uniform estimates in the sequel.
The space $\mH^{1}_{0}(\omega)$ will refer to the closure of 
$\mathscr{D}(\omega):=\{ \varphi\in \mathscr{C}^{\infty}(\RR^{d}), \;
\supp(\varphi)\subset\omega, \; \supp(\varphi)\;\text{bounded}\}$
for $\Vert \cdot \Vert_{\mH^{1}(\omega)}$.

\quad\\
Next we introduce the space of Dirichlet traces
$\mH^{1/2}(\partial\omega):=\{ v\vert_{\partial\omega}, \;
v\in \mH^{1}(\RR^{d})\}$ equipped with the quotient norm
$\Vert v\Vert_{\mH^{1/2}(\partial\omega)}:=\min\{\Vert \varphi\Vert_{\mH^{1}(\RR^{d})},
\varphi\in\mH^{1}(\RR^{d})\,\text{and}\,\varphi\vert_{\partial\omega} = v\}$.
The topological dual to $\mH^{1/2}(\partial\omega)$ will be denoted
$\mH^{-1/2}(\partial\omega) = \mH^{1/2}(\partial\omega)^{*}$.
As detailed for example in \cite[Thm.3.38]{zbMATH01446717}, the trace
map gives rise to a bounded linear operator
\begin{equation}\label{DefinitionTraceOperator}
  \begin{aligned}
    & \mB_{\omega}:\mH^{1}(\omega)\to \mH^{1/2}(\partial\omega)\\
    & \mB_{\omega}(v):= v\vert_{\partial\omega}\quad \forall v\in \mathscr{D}(\RR^{d}).
  \end{aligned}
\end{equation}
We underline that $\mB_{\omega}$ refers to the trace taken from the \textit{interior}
of $\omega$. The norm \eqref{H1fctspace} gives rise to a natural right-inverse of this
Dirichlet boundary trace operator. We define the harmonic lifting operator
$\mB_{\omega}^{\dagger}:\mH^{1/2}(\partial\omega)\to \mH^{1}(\omega)$, see
\cite[\S 1.2.2.4]{MR3013465}, through norm minimization
\begin{equation}\label{HarmonicLifting}
  \begin{aligned}
    & \mB_{\omega}\cdot\mB^{\dagger}_{\omega}(v) = v\quad\forall v\in\mH^{1/2}(\partial\omega)\;\; \text{and}\\
    & \Vert \mB^{\dagger}_{\omega}(v)\Vert_{\mH^{1}(\omega)}:= \min\{
    \Vert \phi\Vert_{\mH^{1}(\omega)},\; \mB_{\omega}(\phi) = v, \;\phi\in \mH^{1}(\omega)\}.\\[2pt]
  \end{aligned}
\end{equation}
Denote $\mH^{1}(\Delta,\omega):=\{v\in\mH^{1}(\Omega),\Delta v\in\mL^{2}(\Omega)\}$ and
let $\bn_{\omega}$ refer to the unit normal vector field to the boundary $\partial\omega$
directed toward the exterior of $\omega$. The Dirichlet trace operator
$\varphi\mapsto \varphi\vert_{\partial\omega}$, resp. the Neumann trace operator
$\varphi\mapsto \bn_{\omega}\cdot\nabla \varphi\vert_{\partial\omega}$, can be extended
by density as a bounded linear map $\mH^{1}(\omega)\to \mH^{1/2}(\partial\omega)$ resp.
$\mH^{1}(\Delta,\omega)\to \mH^{-1/2}(\partial\omega)$, see e.g.
\cite[Lem.4.3]{zbMATH01446717}. The Dirichlet-to-Neumann (DtN) map
$\mT_{\omega}:\mH^{1/2}(\partial\omega)\to \mH^{-1/2}(\partial\omega)$
is defined as the unique bounded linear operator satisfying 
\begin{equation}\label{DefinitionDTN}
  \begin{aligned}
    & \mT_{\omega}(\phi\vert_{\partial\omega}):= \bn_{\omega}\cdot\nabla\phi\vert_{\partial\omega}\\
    & \forall\phi \in\mH^{1}(\Delta,\omega)\;\text{satisfying}\\
    & -\Delta \phi + \gamma^{-2}\phi = 0\quad \text{in}\;\omega.
  \end{aligned}
\end{equation}
This is a real valued and self-adjoint operator $\overline{\mT}_{\omega} =
\mT_{\omega}$ and $\mT_{\omega}^{*} = \mT_{\omega}$ which
induces a scalar product over $\mH^{+1/2}(\partial\omega)$
and the Neumann-to-Dirichlet map $\mT^{-1}_{\omega}: \mH^{-1/2}(\partial\omega)\to
\mH^{+1/2}(\partial\omega)$ induces a scalar product over $\mH^{-1/2}(\partial\omega)$.
We set 
\begin{equation}\label{NormImpedance}
  \begin{aligned}
    & \Vert v\Vert_{\mT_{\omega}}^{2}:=
    \langle \mT_{\omega}(v),\overline{v}\rangle\\
    & \Vert p\Vert_{\mT_{\omega}^{-1}}^{2}:=
    \langle \mT_{\omega}^{-1}(p),\overline{p}\rangle.
  \end{aligned}
\end{equation}
It is a well established fact (see e.g. \cite[Def.1.41]{MR3013465}
or \cite[\S 6.6.3]{MR2361676}) that $\Vert \cdot\Vert_{\mH^{1/2}(\partial\omega)}$
and $\Vert \cdot\Vert_{\mH^{-1/2}(\partial\omega)}$ are  equivalent to the norms \eqref{NormImpedance}.
Applying the Euler equation characterizing the harmonic lifting $\mB^{\dagger}_{\omega}(v)$
as unique solution to the minimization \eqref{HarmonicLifting}, see e.g.
\cite[Thm.7.2-1]{zbMATH00195024}, we have
$-\Delta \mB^{\dagger}_{\omega}(v)+\gamma^{-2}\mB^{\dagger}_{\omega}(v) = 0$ in $\omega$,
so that $\mT_{\omega}(v) = \bn_{\omega}\cdot\nabla \mB^{\dagger}(v)\vert_{\partial\omega}$.
We also deduce that $  \Vert \phi\vert_{\partial\omega}\Vert_{\mT_{\omega}} =
\Vert \mB_{\omega}^{\dagger}(\phi\vert_{\partial\omega})\Vert_{\mH^{1}(\omega)} \leq
\Vert \phi\Vert_{\mH^{1}(\omega)}$ for all $\phi\in\mH^{1}(\omega)$ and, in particular,
we have the inequalities
\begin{equation}\label{ContinuiteOperateurTrace}
  \begin{aligned}
    & \Vert \mB_{\omega}^{\dagger}(v)\Vert_{\mH^1(\omega)} =
    \Vert v\Vert_{\mT_{\omega}}
    \quad \forall v\in\mH^{1/2}(\partial\omega),\\
    & \Vert \mB_{\omega}(u)\Vert_{\mT_{\omega}}\leq \Vert u\Vert_{\mH^{1}(\omega)}
    \quad \forall u\in\mH^{1}(\omega).\\
  \end{aligned}
\end{equation}

\section{Single domain variational formulation}\label{SingleDomainVarForm}

The next step in our analysis will consist in writing Problem \eqref{PbInit1}
in a variational form able to cope with a variety of boundary conditions.
This is why we treat the boundary condition by means of an additional Lagrange
parameter. Let $\Omega\subset \RR^d,\;\;\Gamma:=\partial\Omega$
refer to an open bounded Lipschitz set and its boundary and denote
\begin{equation*}
  \mH(\Omega\times\Gamma):=\mH^{1}(\Omega)\times\mH^{-1/2}(\Gamma)
\end{equation*}
Our analysis will start from a variational formulation
of \eqref{PbInit1}, later referred to as the primary formulation,
that we write: find $\bu\in \mH(\Omega\times\Gamma)$ such that 
\begin{equation}\label{PrimaryFormulation}
    \mA_{\Omega\times\Gamma}(\bu) = \ell_{\Omega\times\Gamma}
\end{equation}
where the bilinear map underlying the variational problem is written as
a bounded linear operator $\mA_{\Omega\times\Gamma}:\mH(\Omega\times\Gamma)\to
\mH(\Omega\times\Gamma)^{*}$ assumed to systematically take
the following form: for any $u,v\in\mH^{1}(\Omega)$ and $p,q\in\mH^{-1/2}(\Gamma)$,
\begin{empheq}[box=\widefbox]{equation}\label{UnifiedVariationalSetting}\tag{\text{A1}}
  \begin{array}{c}
  \textbf{Assumption:}\\[5pt]
     \langle  \mA_{\Omega\times\Gamma}(u,p),(v,q) \rangle:=
     \langle  \mA_{\Omega}(u),v\rangle + \langle \mA_{\Gamma}(u\vert_{\Gamma},p),
     (v\vert_{\Gamma},q)\rangle\\[2pt]
  \end{array}
\end{empheq}
The map $\mA_{\Omega\times\Gamma}$ involves a volume part $\mA_{\Omega}:\mH^{1}(\Omega)\to
\mH^{1}(\Omega)^{*}$ that accounts for the Helmholtz equation in the interior of the domain
$\Omega$. For  $\mu\in\CC$ and $\kappa:\Omega\to \mathbb{C}$
an essentially bounded measurable function, it is assumed of the following form
\begin{empheq}[box=\widefbox]{equation}\label{SuperBilinearForm}\tag{\text{A2}}
  \begin{array}{c}
    \textbf{Assumptions:}\\[5pt]
    \langle  \mA_{\Omega}(u),v\rangle :=
    \int_{\Omega}\mu^{-1}\nabla u\cdot\nabla v - \kappa^{2} uv\,d\bx,\\[5pt]
    \text{with}\;\;\Im m\{\kappa(\bx)^{2}\}\geq 0, \;\forall \bx\in \Omega\\
    \mrm{sup}_{\bx\in \Omega}\vert\kappa(\bx)\vert<\infty\\ 
    \Re e\{\mu\}>0,\;\Im m\{\mu\}\geq 0.
  \end{array}
\end{empheq}
The assumptions above imply in particular that $\Im m\{\langle \mA_{\Omega}(u),\overline{u}\rangle\}
\leq 0\;\forall u\in\mH^{1}(\Omega)$. The operator $\mA_{\Omega\times\Gamma}$ also involves
a pure boundary part $\mA_{\Gamma}$ that
models boundary conditions,
\begin{equation}
  \begin{aligned}
    & \mA_{\Gamma}:\mH^{\textsc{b}}(\Gamma)\to \mH^{\textsc{b}}(\Gamma)^*\\
    & \text{where}\;\mH^{\textsc{b}}(\Gamma):=\mH^{1/2}(\Gamma)\times \mH^{-1/2}(\Gamma).
  \end{aligned}
\end{equation}
The boundary operator $\mA_{\Gamma}$ involves traces on $\Gamma$ and is chosen in
accordance with the boundary conditions of our primary boundary value problem \eqref{PbInit1}.
We will need to rely on the following additional assumptions
\begin{empheq}[box=\widefbox]{equation}\label{AssumptionAbsorption}\tag{\text{A3}}
  \begin{array}{c}
    \textbf{Assumptions:}\\[5pt]
    \textit{i)}\;\;
    \Im m\{\langle \mA_{\Gamma}(\bu),\overline{\bu}\rangle\}\leq 0\quad
    \forall \bu\in \mH^{\textsc{b}}(\Gamma)\\
    \textit{ii)}\;\;
    \mrm{range}(\mA_{\Omega\times\Gamma})\;\text{is closed in }\mH(\Omega\times\Gamma)^*.\\[2pt]
  \end{array}
\end{empheq}
In the remaining of this contribution we will almost systematically take
\eqref{UnifiedVariationalSetting}-\eqref{SuperBilinearForm}-\eqref{AssumptionAbsorption}
as assumptions. We do not require that $\mA_{\Omega\times\Gamma} = \mA_{\Omega\times\Gamma}^{*}$. 
Let us underline that the assumptions above are fulfilled by
$\mA_{\Omega}, \mA_{\Gamma}, \mA_{\Omega\times\Gamma}$ if and only if they are fulfilled
by $\mA_{\Omega}^{*}, \mA_{\Gamma}^{*}, \mA_{\Omega\times\Gamma}^{*}$ (recall that adjunction
does not involve any complex conjugation here). The last hypothesis in
\eqref{AssumptionAbsorption} implies (see e.g. \cite[Thm.2.19]{MR2759829})
\begin{equation}\label{FredholmAlternative}
  \mrm{range}(\mA_{\Omega\times\Gamma}) =
  \mrm{ker}(\mA_{\Omega\times\Gamma}^{*})^{\circ}.
\end{equation}
hence $\mrm{codim}(\mrm{range}(\mA_{\Omega\times\Gamma})) = \mrm{dim}(\mrm{ker}(\mA_{\Omega\times\Gamma}^{*}))$. 
The source functional in \eqref{PrimaryFormulation} is assumed to take the similar
form $\langle \ell_{\Omega\times\Gamma},(v,q)\rangle:=\langle \ell_{\Omega},v\rangle +\langle
\ell_{\Gamma}, (v\vert_{\Gamma},q)\rangle$ where $\langle \ell_{\Omega},v\rangle:=\int_{\Omega}
fv \,d\bx$ for some $f\in\mL^{2}(\Omega)$ and $\ell_{\Gamma}\in \mH^{\textsc{b}}(\Gamma)^* =
\mH^{-1/2}(\Gamma)\times \mH^{+1/2}(\Gamma)$ is chosen in accordance with the boundary condition.

\quad\\
Now we consider concrete boundary conditions,
exhibit corresponding appropriate choices of $\mA_{\Gamma}$
and point how these situations fit the previous assumptions 
\eqref{UnifiedVariationalSetting}-\eqref{SuperBilinearForm}-\eqref{AssumptionAbsorption}.
Here and in the following, for the sake of conciseness, we shall take the notational
convention (see \eqref{DefinitionDTN}),
\begin{equation*}
  \mT_{\Gamma} := \mT_{\RR^{d}\setminus\overline{\Omega}}.
\end{equation*}

\begin{example}[\textbf{Dirichlet boundary condition}]\label{DirichletBC1}
  In the case of a Dirichlet boundary condition, we set
  $\mA_{\Gamma}(\alpha,p) := (p,\alpha)$ and $\ell_{\Gamma}:=(0,g)$
  for some $g\in\mH^{1/2}(\Gamma)$. We have $\Im m\{\langle
  \mA_{\Gamma}(\bu),\overline{\bu}\rangle\} = 0$ for all $\bu\in
  \mH^{\textsc{b}}(\Gamma)$, which fits \textit{i)} of
  \eqref{AssumptionAbsorption}. Formulation \eqref{PrimaryFormulation}
  formulation of a Helmholtz problem with a Dirichlet condition
  imposed by means of a Lagrange parameter at the boundary
  \begin{equation*}
    \begin{array}{ll}
      u\in \mH^{1}(\Omega),\;p\in \mH^{-1/2}(\Gamma)\;\text{such that}\\
      \int_{\Omega}\mu^{-1}\nabla u\cdot\nabla v - \kappa^{2} uv\,d\bx +
      \int_{\Gamma}p v\, d\sigma = \int_{\Omega}fv d\bx & \forall v\in \mH^{1}(\Omega),\\
      \int_{\Gamma}u q \, d\sigma = \int_{\Gamma}g q\,d\sigma & \forall q\in \mH^{-1/2}(\Gamma).   
    \end{array}
  \end{equation*}
  Whenever there is existence and uniqueness of the solution pair $(u,p)$ then
  $p = -\bn_{\Omega}\cdot\nabla u \vert_{\Gamma}$. Conditions in \eqref{SuperBilinearForm}
  guarantee that the volume part of this equation is coercive modulo the compact term
  attached to $\kappa$. Hence the operator associated to this
  system is of Fredholm type with index $0$. In particular it has closed range,
  which fits \textit{ii)} of \eqref{AssumptionAbsorption}.
\end{example}

\begin{example}[\textbf{Neumann boundary condition}]\label{NeumannBC1}
  In the case of Neumann conditions, the boundary data is
  $g\in\mH^{-1/2}(\Gamma)$ and we choose $\mA_{\Gamma} (\alpha,p) :=
  (0,\mT_{\Gamma}^{-1}p)$ and $\ell_{\Gamma}:= (g,0)$. Again we have
  $\Im m\{\langle \mA_{\Gamma}(\bu),\overline{\bu}\rangle\} = 0$ for
  all $\bu\in\mH^{\textsc{b}}(\Gamma)$, so this choice also matches
  \textit{i)} of \eqref{AssumptionAbsorption}.  The primary
  formulation \eqref{PrimaryFormulation} writes
  \begin{equation}\label{NeumannFormulation}
    \begin{array}{ll}
      u\in \mH^{1}(\Omega),\;p\in \mH^{-1/2}(\Gamma)\;\text{such that}\\
      \int_{\Omega}\mu^{-1}\nabla u\cdot\nabla v - \kappa^{2} uv\,d\bx
      = \int_{\Omega} fv d\bx + \int_{\Gamma} g vd\sigma & \forall v\in \mH^{1}(\Omega),\\
      \int_{\Gamma}q\mT_{\Gamma}^{-1} p \, d\sigma = 0    & \forall q\in \mH^{-1/2}(\Gamma),
    \end{array}
  \end{equation}
  where $u$ is decoupled from $p$. Actually we have in particular
  $p=0$ and this variable is not supposed to receive any particular
  interpretation.  Since $\mT_{\Gamma}^{-1}:\mH^{-1/2}(\Gamma)\to
  \mH^{1/2}(\Gamma)$ is an isomorphism, the operator
  $\mA_{\Omega\times\Gamma}$ associated to \eqref{NeumannFormulation}
  is of Fredholm type with index 0.
\end{example}

\begin{example}[\textbf{Robin boundary condition}]\label{RobinBC1}
Consider a bounded linear map $\Lambda:\mH^{1/2}(\Gamma)\to \mH^{-1/2}(\Gamma)$
that satisfies $\Re e\{\langle \Lambda(v),\overline{v}\rangle\} >0
\;\forall v\in \mH^{1/2}(\Gamma)\setminus\{0\}$ (as a typical example:
$\Lambda(v) = \lambda v$ with $\lambda>0$). 
In this case again the boundary data is $g\in\mH^{-1/2}(\Gamma)$ and we choose
$\mA_{\Gamma} (\alpha,p) :=(-i\Lambda\alpha,\mT_{\Gamma}^{-1}p)$ and
$\ell_{\Gamma}:=(g,0)$.
This choice of $\mA_{\Gamma}$ corresponds to the boundary condition
$\bn_{\Omega}\cdot\nabla u\vert_{\Gamma}-i\Lambda(u) = 0$ on $\Gamma$. 
Formulation \eqref{PrimaryFormulation} writes
\begin{equation*}
  \begin{array}{ll}
    u\in \mH^{1}(\Omega),\;p\in \mH^{-1/2}(\Gamma)\;\text{such that}\\
    \int_{\Omega}\mu^{-1}\nabla u\cdot\nabla v - \kappa^{2} uv\,d\bx - i\int_{\Gamma}v\Lambda(u) d\sigma
    = \int_{\Omega} fv d\bx + \int_{\Gamma} g vd\sigma & \forall v\in \mH^{1}(\Omega)\\
    \int_{\Gamma}q\mT_{\Gamma}^{-1} p \, d\sigma = 0    & \forall q\in \mH^{-1/2}(\Gamma)
  \end{array}
\end{equation*}
which is a variant of \eqref{NeumannFormulation} involving 
$i\int_{\Gamma}v\Lambda(u) d\sigma$ as an additional term.
Again $p$ is decoupled from the rest of the system and $p=0$.
Again the operator $\mA_{\Omega\times\Gamma}$ associated to this system is of Fredholm type
with index 0.
\end{example}

\begin{example}[\textbf{Mixed boundary conditions}]\label{MixedBC1}
Assume that the boundary $\Gamma$ is decomposed in two Lipschitz
  regular parts $\Gamma = \overline{\Gamma}_{\dir}\cup
  \overline{\Gamma}_{\neu}$ where $\Gamma_{\dir},\Gamma_{\neu}$ have
  non-vanishing surface measure. Suppose we wish to consider mixed
  Dirichlet-Neumann conditions
  \begin{equation*}
    u\vert_{\Gamma_{\dir}} = g_{\dir}\vert_{\Gamma_{\dir}}\quad\text{and}\quad
    \bn_{\Omega}\cdot\nabla u\vert_{\Gamma_{\neu}} = g_{\neu}\vert_{\Gamma_{\neu}}
  \end{equation*}
  for some $g_{\dir} \in \mH^{1/2}(\Gamma)$ and $g_{\neu} \in \mH^{-1/2}(\Gamma)$. Defining
  $\widetilde{\mH}^{-1/2}(\Gamma_{\dir}):=\{q\in \mH^{-1/2}(\Gamma),\;q\vert_{\Gamma_{\neu}} = 0\}$,
  the Helmholtz problem with the mixed boundary conditions above is then equivalent to the
  variational formulation
  \begin{equation*}
    \begin{array}{ll}
      u\in \mH^{1}(\Omega),\;\tilde{p}\in \widetilde{\mH}^{-1/2}(\Gamma_{\dir})\;\text{such that}\\
      \int_{\Omega}\mu^{-1}\nabla u\cdot\nabla v - \kappa^{2} uv\,d\bx + \int_{\Gamma}\tilde{p} v d\sigma
      = \int_{\Omega} fv d\bx + \int_{\Gamma} g_{\neu} vd\sigma & \forall v\in \mH^{1}(\Omega)\\
      \int_{\Gamma}u \tilde{q} d\sigma = \int_{\Gamma}g_{\dir} \tilde{q} d\sigma
      & \forall \tilde{q}\in \widetilde{\mH}^{-1/2}(\Gamma_{\dir}).
    \end{array}    
  \end{equation*}
  The space $\widetilde{\mH}^{-1/2}(\Gamma_{\dir})$ is closed in
  $\mH^{-1/2}(\Gamma)$.  Considering the
  $\mT_{\Gamma}^{-1}$-orthogonal projector
  $\Theta:\mH^{-1/2}(\Gamma)\to \mH^{-1/2}(\Gamma)$ with
  $\mrm{range}(\Theta) = \widetilde{\mH}^{-1/2}(\Gamma_{\dir})$, the
  previous variational formulation rewrites equivalently as
  \begin{equation*}
    \begin{array}{ll}
      u\in \mH^{1}(\Omega),\;p\in \mH^{-1/2}(\Gamma)\;\text{such that}\\
      \int_{\Omega}\mu^{-1}\nabla u\cdot\nabla v - \kappa^{2} uv\,d\bx + \int_{\Gamma}v \Theta p\, d\sigma
      = \int_{\Omega} fv d\bx + \int_{\Gamma} g_{\neu} vd\sigma & \forall v\in \mH^{1}(\Omega)\\
      \int_{\Gamma}q\Theta^{*}u d\sigma + \int_{\Gamma}q\,\mT_{\Gamma}^{-1}(\Id-\Theta)\,p d\sigma
      = \int_{\Gamma}q\Theta^{*}g_{\dir} d\sigma
      & \forall q\in \mH^{-1/2}(\Gamma).
    \end{array}    
  \end{equation*}
  This shows that the case of mixed Dirichlet-Neumann boundary conditions fits into
  the framework \eqref{UnifiedVariationalSetting}-\eqref{SuperBilinearForm} with the
  choice $\mA_{\Gamma} (\alpha,p) :=(\Theta p, \Theta^{*}\alpha + \mT_{\Gamma}^{-1}(\Id-\Theta) p)$ and
  $\ell_{\Gamma}:=(g_{\neu},\Theta^{*}g_{\dir})$. The $\mT_{\Gamma}^{-1}$-orthogonality of $\Theta$
  yields $\mT_{\Gamma}^{-1}(\Id-\Theta) = (\Id-\Theta^{*})\mT_{\Gamma}^{-1}(\Id-\Theta)$
  and $\Theta(\overline{p}) = \overline{\Theta(p)}$, so that 
  $\langle \mA_{\Gamma} (\alpha,p), (\overline{\alpha},\overline{p})\rangle =
  \langle\overline{\alpha},\Theta(p)\rangle + \langle\overline{\Theta(p),\overline{\alpha}}\rangle +
  \langle \overline{(\Id-\Theta)p},\mT_{\Gamma}^{-1}(\Id-\Theta) p\rangle $.
  We conclude that $\Im m\{\langle \mA_{\Gamma} (\bu), \overline{\bu}\rangle\} = 0$ for all
  $\bu\in \mH^{\textsc{b}}(\Gamma)$. 
\end{example}

\section{Multi-domain functional setting}
The boundary value problem \eqref{PbInit1} has been reformulated as an
equivalent global variational problem with \eqref{PrimaryFormulation}.
As we aim at extending an analytical framework for domain decomposition
by substructuration though, we are going to reshape Formulation \eqref{PrimaryFormulation},
adapting it to a multi-domain geometrical configuration. For this,
we need to introduce notations adapted to domain decomposition.
Consider a decomposition into a collection of non-overlapping Lipschitz
open sets $\Omega_{j}\subset \RR^{d},j=1,\dots,\mJ$ that satisfy
\begin{equation}\label{SubdomainPartition}
  \begin{aligned}
    & \overline{\Omega} = \overline{\Omega}_1 \cup \dots\cup \overline{\Omega}_\mJ,\\
    & \text{with}\;\Omega_j\cap \Omega_k = \emptyset\;\text{for}\;j\neq k.
  \end{aligned}
\end{equation}
Such a decomposition may very well admit a non-trivial wire-basket
i.e. the set of cross points is non-empty, and we wish to underline that
this situation is covered by the subsequent analysis. We shall
refer to the skeleton of the decomposition by 
\begin{equation}
  \Sigma:= \partial\Omega_1\cup\dots\cup\partial\Omega_\mJ.
\end{equation}
Note that $\Gamma = \partial\Omega\subset \Sigma$. We need to introduce notations for function
spaces adapted to this multi-domain setting. In this context, cartesian product spaces are probably
the most natural, so we set 
\begin{equation}\label{MultiDomainSpaces}
  \begin{aligned}
    \mH^{\textsc{b}}(\Gamma) &:= \mH^{\frac{1}{2}}(\Gamma)\times \mH^{-\frac{1}{2}}(\Gamma)\\
    \mbH(\Omega) & := \mH^{\textsc{b}}(\Gamma)
    \times \mH^{1}(\Omega_1)\times\dots\times  \mH^{1}(\Omega_\mJ)\\
    \mbH(\Sigma) & :=\mH^{\frac{1}{2}}(\Gamma)\times \mH^{\frac{1}{2}}(\partial\Omega_{1})
    \times \dots\times \mH^{\frac{1}{2}}(\partial\Omega_{\mJ})
  \end{aligned}
\end{equation}
As cartesian products, these spaces are equipped with norms and duality pairings
given by \eqref{CartesianProductNorm}. Apart from the boundary terms attached to
$\mH^{\textsc{b}}(\Gamma)$,  the space $\mbH(\Omega)$ should be understood as functions
defined over $\Omega$, admitting potential jumps through interfaces.
The space $\mbH(\Sigma)$ consists in tuples of Dirichlet traces. Its dual is 
\begin{equation*}
  \mbH(\Sigma)^* = \mH^{-\frac{1}{2}}(\Gamma)\times  \mH^{-\frac{1}{2}}(\partial\Omega_1) \times \dots\times
  \mH^{-\frac{1}{2}}(\partial\Omega_\mJ).
\end{equation*}
We need to introduce several operators acting in these spaces. First we shall consider
the operator $\mT:\mbH(\Sigma)\to \mbH(\Sigma)^{*}$ defined as the block
diagonal operator acting locally in each subdomain
\begin{equation}\label{DefGlobalDtN}
  \begin{aligned}
    & \mT:= \mrm{diag}(\mT_{\Gamma},\mT_{\Omega_{1}},\dots,\mT_{\Omega_{\mJ}})\\
    & \text{where}\;\;\mT_{\Gamma}:=\mT_{\RR\setminus\overline{\Omega}}
  \end{aligned}
\end{equation}
and each $\mT_{\Omega_j}$ is defined with \eqref{DefinitionDTN}.
The norms $\Vert \cdot\Vert_{\mT}$ and  $\Vert \cdot\Vert_{\mT^{-1}}$
defined by \eqref{OperatorInducedNorm} and \eqref{DefGlobalDtN} are equivalent to
$\Vert \cdot\Vert_{\mbH(\Sigma)}$ and  $\Vert \cdot\Vert_{\mbH(\Sigma)^{*}}$,
which stems from the analogous property being satisfied locally by each
$\mT_{\Omega_j}$. These norms will play an important role in the subsequent analysis.
Next we introduce a boundary trace operator
$\mB:\mbH(\Omega)\to \mbH(\Sigma)$ and defined by
\begin{equation}
  \begin{aligned}
    & \mB := \mrm{diag}(\mB_{\Gamma},\mB_{\Omega_1},\dots,\mB_{\Omega_\mJ})\\
    & \text{where}\;\; \mB_{\Gamma}(\alpha,p):=\alpha
  \end{aligned}
\end{equation}
and each $\mB_{\Omega_j}$ is the Dirichlet trace operator interior to subdomain $\Omega_j$
as defined in \eqref{DefinitionTraceOperator}. By definition of $\mT$ we have
$\Vert\mB(\bu)\Vert_{\mT}\leq \Vert \bu\Vert_{\mbH(\Omega)}$ for all
$\bu\in\mbH(\Omega)$, since a similar inequality is satisfied in each subdomain
locally according to \eqref{ContinuiteOperateurTrace}. We can also form a
multi-domain harmonic lifting map $\mB^{\dagger}:\mbH(\Sigma)\to \mbH(\Omega)$
defined as the block-diagonal operator as follows
\begin{equation}
  \begin{aligned}
    & \mB^{\dagger} = \mrm{diag}(\mB_{\Gamma}^{\dagger}, \mB_{\Omega_1}^{\dagger},
    \dots,\mB_{\Omega_\mJ}^{\dagger})\\
    & \text{where}\;\; \mB_{\Gamma}^{\dagger}(\alpha):= (\alpha,0)\\
  \end{aligned}
\end{equation}
and each $\mB_{\Omega_j}^{\dagger}$ as defined in \eqref{HarmonicLifting}.
With this definition we have $\mB\mB^{\dagger} = \mrm{Id}$ and $\mB^{\dagger}\mB$
is an orthogonal projector in $\mbH(\Omega)$. Finally we also need to consider a restriction
operator $\mR:\mH(\Omega\times\Gamma)\to \mbH(\Omega)$ that embeds pairs
$(u,p)\in \mH(\Omega\times\Gamma) = \mH^{1}(\Omega)\times \mH^{-1/2}(\Gamma)$
into the cartesian product $\mbH(\Omega)$ by restricting locally to each
subdomain
\begin{equation}\label{DefOperatorRestriction}
  \begin{aligned}
    & \mR(u,p) := ((u\vert_{\Gamma},p),u\vert_{\Omega_1},\dots,u\vert_{\Omega_{\mJ}})\\
    & \text{for}\; u\in \mH^{1}(\Omega), p\in\mH^{-1/2}(\Gamma).
  \end{aligned}
  \end{equation}
The image of this operator $\mrm{range}(\mR) = \mR(\mH(\Omega\times\Gamma))$ is a
particular subspace of $\mbH(\Omega)$ spanned by tuples of functions
that match through interfaces. This matching property is precisely what
characterizes Dirichlet transmission conditions through
interfaces of the decomposition \eqref{SubdomainPartition}. This is why we
dedicate notations to this. 
\begin{equation}\label{DefinitionSingleTraceSpace}
  \begin{aligned}
    & \mbX(\Omega):=\{ \mR(u,p),\;u\in \mH^{1}(\Omega),p\in \mH^{-1/2}(\Gamma)\}\\
    & \mbX(\Sigma):=\{ \mB(\bu),\;\bu\in \mbX(\Omega)\}\\
    & \mbX(\Sigma)^{\circ}:=\{\bp\in \mbH(\Sigma)^{*},\;\langle \bp,\bv\rangle = 0\;\forall \bv\in\mbX(\Sigma)\}.
  \end{aligned}
\end{equation}
A rapid inspection of the previous definitions shows that $\mbX(\Sigma) = \{ (u\vert_{\Gamma},
u\vert_{\partial\Omega_1},\dots,u\vert_{\partial\Omega_{\mJ}}), u\in \mH^{1}(\Omega)\}$ i.e. these are the tuples
of Dirichlet traces that match through interfaces. The space $\mbX(\Sigma)$ (resp. $\mbX(\Omega)$)
is a closed subspace of $\mbH(\Sigma)$ (resp. $\mbH(\Omega)$) that encodes the Dirichlet
transmission conditions through interfaces, while $\mbX(\Sigma)^{\circ}$ is a closed
subspace of $\mbH(\Omega)^{*}$ that encodes the Neumann transmission conditions.
Indeed, considering restriction to interfaces in the sense of distributions, 
\begin{equation}
  \begin{aligned}
    & (v_0,\dots, v_{\mJ})\in \mbX(\Sigma)^{\textcolor{white}{\circ}}\;
    \Longrightarrow\; v_j = +v_k\;\text{on}\;\Gamma_j\cap\Gamma_k,\\
    & (p_0,\dots, p_{\mJ})\in \mbX(\Sigma)^{\circ}\;
    \Longrightarrow\; p_j = -p_k\;\text{on}\;\Gamma_j\cap\Gamma_k.
  \end{aligned}
\end{equation}
It is clear from these definitions that $\mbX(\Omega) = \{\bu\in \mbH(\Omega),
\;\mB(\bu)\in\mbX(\Sigma)\}$. In particular $\mrm{ker}(\mB)\subset \mbX(\Omega)$.
Recall the definition of polar sets given by \eqref{DefPolarSet}.
The following lemma is a continuous counterpart to \cite[Lem.2.1]{claeys2021nonself}.
\newpage 

\begin{lem}\label{ExtendedPolarity}
  \quad
  \begin{itemize}
  \item[i)]    $\mrm{ker}(\mB)^{\circ} = \mrm{range}(\mB^*)$\\[-20pt]
  \item[ii)]   $\mrm{ker}(\mB^{*}) = \{0\}$\\[-20pt]
  \item[iii)]  $\mbX(\Omega) =\mB^{-1}(\mbX(\Sigma))$\\[-20pt]
  \item[iv)]   $\mbX(\Omega)^{\circ} = \mB^*(\mbX(\Sigma)^{\circ})$
  \end{itemize}
\end{lem}
\noindent \textbf{Proof:}

The first and second results are direct consequences of the surjectivity of the trace map
$\mB:\mbH(\Omega)\to \mbH(\Sigma)$ combined with Theorem 4.7, 4.12 and 4.15
of \cite{zbMATH01022519}. The third result is a rephrasing of
$\mbX(\Omega) = \{\bu\in \mbH(\Omega),\;\mB(\bu)\in\mbX(\Sigma)\}$ in condensed form.
To prove the last result, first observe that $\mB^*(\mbX(\Sigma)^{\circ})\subset
\mbX(\Omega)^{\circ}$ by routine verifications.

Now pick an arbitrary $\bp\in \mbX(\Omega)^{\circ}$.
Since $\mrm{ker}(\mB)\subset \mbX(\Omega)\Rightarrow \mbX(\Omega)^{\circ}
\subset \mrm{ker}(\mB)^{\circ} = \mrm{range}(\mB^*)$, there exists
$\bq\in \mbH(\Sigma)^{*}$ such that $\bp = \mB^{*}\bq$. For any $\bv\in \mbX(\Sigma)$,
there exists $\bu\in \mbX(\Omega)$ such that $\bv = \mB(\bu)$,
which implies that $\langle \bq,\bv\rangle =\langle \bp,\bu\rangle = 0$. 
From this we conclude that $\bq\in \mbX(\Sigma)^{\circ}$ hence $\bp\in \mB^{*}(\mbX(\Sigma)^{\circ})$,
which proves $\mbX(\Omega)^{\circ}\subset \mB^{*}(\mbX(\Sigma)^{\circ})$. \hfill $\Box$

\quad\\
In Item \textit{iii)} of the lemma above, $\mB^{-1}(\mbX(\Sigma)) = \{\bu\in\mbH(\Omega),
\mB(\bu)\in \mbX(\Sigma)\}$ refers to a pre-image (the operator $\mB$ is obviously
non-invertible i.e. $\mrm{ker}(\mB)\neq \{0\}$). The following orthogonal decomposition
was established in \cite[Prop.4.2]{zbMATH01446717}.

\begin{prop}\label{OrthogonalSum}\quad\\
  We have $\mbH(\Sigma)^* = \mbX(\Sigma)^{\circ}\oplus \mT(\mbX(\Sigma))$
  and this decomposition is $\mT^{-1}$-orthogonal.
\end{prop}

\noindent 
The orthogonal decomposition of the previous result can be used
to elaborate a characterization of transmission conditions.
The following result was established in \cite[Prop.5.4]{zbMATH01446717}. 

\begin{prop}\label{NonlocalTrCond}\quad\\
  Let $\mQ:\mbH(\Sigma)^{*}\to \mbH(\Sigma)^{*}$ refer to the $\mT^{-1}$-orthogonal projection
  onto  $\mT(\mbX(\Sigma))$. Then the operator $\Pi:=2\mQ -\Id$ is a $\mT^{-1}$-isometric
  involution i.e. $\Pi^{2} = \Id$, $\Vert \Pi(\bq)\Vert_{\mT^{-1}} = \Vert \bq\Vert_{\mT^{-1}}$
  for all $\bq\in\mbH(\Sigma)^*$. Moreover, for any pair
  $(\bu,\bp)\in \mbH(\Sigma)\times\mbH(\Sigma)^*$, we have
  \begin{equation}\label{CaracSingleTraceSpace}
    (\bu,\bp)\in\mbX(\Sigma)\times \mbX(\Sigma)^{\circ}
    \quad\iff\quad -\bp + i\mT(\bu) = \Pi(\bp+i\mT(\bu)).
  \end{equation}
\end{prop}

\noindent 
The characterization above relies on an exchange operator $\Pi$ which
is characteristic of Optimized Schwarz Methods (OSM, see e.g. \cite[Eq.37]{Bendali2006})
and ultra-weak variational formulations (UWVF) see e.g. \cite[Eq.1.19]{MR1618464}.
An explicit expression of this operator in terms of double layer potentials
attached to the operator $-\Delta +\gamma^{-2}$ was provided in \cite[\S 5.2]{claeys2019new}.

\section{Multi-domain variational formulation}

Using the notations introduced in the previous sections, we now
rewrite the primary formulation \eqref{PrimaryFormulation}, decomposing
it according to the subdomain partition \eqref{SubdomainPartition}. Pick $u,v$
arbitrarily in $\mH^{1}(\Omega)$ and expand the integral coming into play in
the definition \eqref{SuperBilinearForm} of $\mA_{\Omega}$. This leads to
\begin{equation}
  \begin{aligned}
    & \langle \mA_{\Omega}u,v\rangle = \langle \mA_{\Omega_1}(u\vert_{\Omega_1}),
    v\vert_{\Omega_1}\rangle + \dots + \langle \mA_{\Omega_{\mJ}}(u\vert_{\Omega_{\mJ}}),
    v\vert_{\Omega_{\mJ}}\rangle\\
    & \text{with}\quad \langle\mA_{\Omega_j}u,v\rangle:=
    \int_{\Omega_j}\mu^{-1}\nabla u \cdot\nabla v - \kappa^{2} uv \,d\bx
  \end{aligned}
\end{equation}
In the expression above only $u\vert_{\Omega_j},v\vert_{\Omega_j}\in \mH^{1}(\Omega_j)$ come into play
in the term attached to $\Omega_j$. The source term in \eqref{PrimaryFormulation} can be
decomposed in a similar manner  $\ell_{\Omega}(v) = \ell_{\Omega_1}(v\vert_{\Omega_1})+\dots
\ell_{\Omega_\mJ}(v\vert_{\Omega_{\mJ}})$. The above decompositions lead to introducing a block-diagonal operator
$\mA:\mbH(\Omega)\to \mbH(\Omega)^*$ associated to these local bilinear forms i.e. defined
by
\begin{equation}\label{DefinitionOperatorMultiDomain}
  \begin{aligned}
    & \mA := \mrm{diag}(\mA_{\Gamma},\mA_{\Omega_1},\dots,\mA_{\Omega_{\mJ}})\\
    & \text{so that}\; \mA_{\Omega\times\Gamma} = \mR^{*}\mA\mR.
  \end{aligned}  
\end{equation}
We have factorized the operator of our primary boundary value problem  $\mA_{\Omega\times\Gamma}$,
and this factorization is interesting from the perspective of domain decomposition because local
subproblems are disconnected from one another in $\mA$. The following property is inherited
from the assumptions we made in \S \ref{SingleDomainVarForm} about $\mA_{\Omega\times \Gamma},\mu,\kappa$
and $\mA_{\Gamma}$,
\begin{equation}\label{GlobalAbsorptionProperty}
  \Im m\{\langle \mA(\bu),\overline{\bu}\rangle\}\leq 0\quad \forall \bu\in \mbH(\Sigma).
\end{equation}
We also need a unique solvability property for local problems with impedance boundary condition.
Because we do not make much specific assumptions regarding the boundary operator $\mA_{\Gamma}$,
we take this further property as an assumption:
\begin{empheq}[box=\widefbox]{equation}\label{LocalUniqueSolvability}\tag{\text{A4}}
  \begin{array}{c}
    \textbf{Assumption:}\\[5pt]
    \mA - i\mB^*\mT\mB:\mbH(\Omega)\to \mbH(\Omega)^*\\
    \text{is an isomorphism.}
  \end{array}
\end{empheq}
A notable consequence of \eqref{SuperBilinearForm}, \eqref{AssumptionAbsorption}
and \eqref{LocalUniqueSolvability} is that $\mrm{ker}(\mA)\cap \mrm{ker}(\mB) = \{0\}$.
Since $\mA,\mT$ and $\mB$ are subdomain-wise block-diagonal, the assumption above
is actually equivalent to imposing that each $\mA_{\Omega_j} - i\mB_{\Omega_j}^*\mT_{\Omega_j}\mB_{\Omega_j}:
\mH(\Omega_j)\to \mH(\Omega_j)^*$ and $\mA_{\Gamma} - i\mB_{\Gamma}^*\mT_{\Gamma}\mB_{\Gamma}:
\mH^{\textsc{b}}(\Gamma)\to \mH^{\textsc{b}}(\Gamma)^*$ are isomorphisms. These conditions 
are fulfilled in many concrete circumstances. As regards interior contributions, for example,
we have the following simple consequence of the unique continuation principle.

\begin{lem}\label{MultiDomainAbsorptionProp}\quad\\
  Assume \eqref{UnifiedVariationalSetting}-\eqref{SuperBilinearForm} and that
  $\mu,\kappa$ are constants (i.e. do not depend on $\bx$).
  Then for any $j = 1,\dots,\mJ$ the operator $\mA_{\Omega_j}-i\mB_{\Omega_j}^*\mT_{\Omega_j}\mB_{\Omega_j}:
  \mH(\Omega_j)\to \mH(\Omega_j)^*$ is an isomorphism.
\end{lem}
\noindent \textbf{Proof:}

Let us denote $\omega = \Omega_j$ for the sake of conciseness.
According to \eqref{SuperBilinearForm}, there exists $\alpha>0$ such that 
\begin{equation*}
  \begin{array}{l}
    \alpha \Vert u\Vert_{\mH^{1}(\omega)}^{2} \leq
    \Re e\{\langle \tilde{\mA}_{\omega}(u),\overline{u}\rangle\}\quad \forall u\in\mH^{1}(\omega),\\
    \langle \tilde{\mA}_{\omega}(u),v\rangle :=
    \langle (\mA_{\omega} - i\mB_{\omega}^*\mT_{\omega}\mB_{\omega})u,v\rangle
    + \int_{\omega}(1+\kappa^{2}) uv d\bx.
  \end{array}
\end{equation*}
Applying Lax-Milgram's lemma, we see that the operator $\tilde{\mA}_{\omega}:\mbH(\omega)\to
\mbH(\omega)^*$ is an isomorphism hence, since it differs by a compact perturbation, that 
$\mA_{\omega} - i\mB_{\omega}^*\mT_{\omega}\mB_{\omega}$ is of Fredholm type with index 0,
see e.g. \cite[Chap.2]{zbMATH01446717}. There only remains to prove that
$\mrm{ker}(\mA_{\omega} - i\mB_{\omega}^*\mT_{\omega}\mB_{\omega}) = \{0\}$.
Pick any $u\in \mH^{1}(\omega)$  such that $(\mA_{\omega}- i\mB_{\omega}^*\mT_{\omega}\mB_{\omega})u = 0$. Then we have
\begin{equation*}
  \Vert \mB_{\omega}(u)\Vert_{\mT_{\omega}}^{2} \leq -\Im m\{\langle(\mA_{\omega}
  - i\mB_{\omega}^*\mT_{\omega}\mB_{\omega})u, \overline{u} \rangle \} = 0.
\end{equation*}
From this we conclude that $u\vert_{\partial\omega} = \mB_{\omega}(u) = 0$ hence 
$\mA_{\omega}(u) = 0$. On the other hand $\mA_{\omega}(u) = 0\Rightarrow
\bn_{\omega}\cdot\nabla u\vert_{\partial\omega} = 0$. There only remains to apply
the unique continuation principle, see e.g. Lemma 2.2 in \cite{zbMATH04136456},
to conclude that $u = 0$ in $\omega$. \hfill $\Box$

\quad\\
Regarding classical boundary conditions and the associated choice
of $\mA_{\Gamma}$, we can also examine the invertibility of
$\mA_{\Gamma}-i\mB_{\Gamma}^*\mT_{\Gamma}\mB_{\Gamma}$.

\begin{example}[\textbf{Dirichlet condition}]\label{DirichletBC2}
  Taking the same notations as in Example \ref{DirichletBC1},
  in this situation we have the following expression
  $(\mA_{\Gamma}-i\mB_{\Gamma}^*\mT_{\Gamma}\mB_{\Gamma})(\alpha,p) =
  (p-i\mT_{\Gamma}\alpha, \alpha)$. We conclude that
  $\mA_{\Gamma}-i\mB_{\Gamma}^*\mT_{\Gamma}\mB_{\Gamma}$ is continuously
  invertible with
  \begin{equation*}
    (\mA_{\Gamma}-i\mB_{\Gamma}^*\mT_{\Gamma}\mB_{\Gamma})^{-1}(p,\alpha) = (\alpha,p+i\mT_{\Gamma}\alpha).
  \end{equation*}
\end{example}

\begin{example}[\textbf{Neumann condition}]\label{NeumannBC2}
  Taking the same notations as in Example \ref{NeumannBC1}, we have
  $(\mA_{\Gamma}-i\mB_{\Gamma}^*\mT_{\Gamma}\mB_{\Gamma})(\alpha,p) =
  (-i\mT_{\Gamma}\alpha, \mT_{\Gamma}^{-1}p)$.  We conclude that
  $\mA_{\Gamma}-i\mB_{\Gamma}^*\mT_{\Gamma}\mB_{\Gamma}$ is continuously invertible with
  \begin{equation*}
    (\mA_{\Gamma}-i\mB_{\Gamma}^*\mT_{\Gamma}\mB_{\Gamma})^{-1}(p,\alpha) = (i\mT_{\Gamma}^{-1}p, \mT_{\Gamma}\alpha).
  \end{equation*}
\end{example}

\begin{example}[\textbf{Robin condition}]\label{RobinBC2}
  Taking the same notations as in Example \ref{RobinBC1}, we have
  $(\mA_{\Gamma}-i\mB_{\Gamma}^*\mT_{\Gamma}\mB_{\Gamma})(\alpha,p) =   (-i(\Lambda + \mT_{\Gamma})\alpha, \mT_{\Gamma}^{-1}p)$.
  Because $\Re e\{\langle \Lambda(\alpha),\overline{\alpha}\rangle\}> 0$
  for all $\alpha\in \mH^{1/2}(\Gamma)$, we see that $\Lambda+\mT_{\Gamma}$ is coercive
  hence invertible and $\mA_{\Gamma}-i\mB_{\Gamma}^*\mT_{\Gamma}\mB_{\Gamma}$ is then continuously
  invertible with
  \begin{equation*}
    (\mA_{\Gamma}-i\mB_{\Gamma}^*\mT_{\Gamma}\mB_{\Gamma})^{-1}(p,\alpha) =
    (i(\Lambda + \mT_{\Gamma})^{-1}p, \mT_{\Gamma}\alpha).
  \end{equation*}
\end{example}

\begin{example}[\textbf{Mixed conditions}]\label{MixedBC2}
Taking the same notations as in Example \ref{MixedBC1}, we have $
  (\mA_{\Gamma}-i\mB_{\Gamma}^*\mT_{\Gamma}\mB_{\Gamma})(\alpha,p) =
  (\Theta p-i\mT_{\Gamma}\alpha, \Theta^{*}\alpha +
  \mT_{\Gamma}^{-1}(\Id-\Theta)p)$.  Recall that $\Theta$ is a
  $\mT_{\Gamma}^{-1}$-ortogonal projector so that
  $\mT_{\Gamma}^{-1}\Theta = \Theta^{*}\mT_{\Gamma}^{-1}=
  \Theta^{*}\mT_{\Gamma}^{-1}\Theta$.  Invertibility of
  $\mA_{\Gamma}-i\mB_{\Gamma}^*\mT_{\Gamma}\mB_{\Gamma}$ is
  then established by readily checking that the inverse is given by the
  following explicit formula
  \begin{equation*}
    \begin{aligned}
      & (\mA_{\Gamma}-i\mB_{\Gamma}^*\mT_{\Gamma}\mB_{\Gamma})^{-1}(p,\alpha) =\\
      & (\Theta^{*}\alpha+i\mT_{\Gamma}^{-1}(\Id-\Theta)p, (\Id-\Theta + i\Theta)\mT_{\Gamma}\alpha + \Theta p).
    \end{aligned}
  \end{equation*}
\end{example}

\quad\\
Similarly to what precedes, define $\bell\in \mbH(\Omega)^*$
by $\langle \bell,\bv\rangle = \ell_{\Gamma}(v_0,q) + \ell_{\Omega_1}(v_1)
+\dots +\ell_{\Omega_\mJ}(v_{\mJ})$, and we have $\ell_{\Omega\times\Gamma} = \mR^{*}\bell$.
The primary variational problem \eqref{PrimaryFormulation}
can then rewritten by means of $\mA$ as follows: find $\bu\in\mH(\Omega\times\Gamma)$
such that $\langle \mA\mR(\bu),\mR(\bv)\rangle = \langle \bell,\mR(\bv)\rangle$ for
all $\bv\in \mH(\Omega\times\Gamma)$. Making use of the definition of $\mbX(\Omega)$
as the image of $\mR$ see \eqref{DefinitionSingleTraceSpace}, this also rewrites 
\begin{equation}\label{DomainDecomposedPb1}
  \begin{aligned}
    & \bu\in \mbX(\Omega)\;\text{and}\\
    & \langle \mA(\bu),\bv\rangle = \langle \bell,\bv\rangle
    \;\; \forall \bv\in \mbX(\Omega).
  \end{aligned}
\end{equation}

\section{Closed linear manifolds interpretation}

Formulation \eqref{PrimaryFormulation} which is the starting point
of this study, is not assumed to be a priori uniquely solvable. The kernel
of $\mA_{\Omega\times\Gamma}$ might be non-trivial. In many relevant applications
though, it is of Fredholm type, and this is why we are interested
in studying how this Fredholmness carries over in the multi-domain context.
For this we are going to consider the skew-symmetric bilinear form
$\lbrack\cdot,\cdot\rbrack:(\;\mbH(\Sigma) \times\mbH(\Sigma)^{*})^{2}\to \CC$
defined by 
\begin{equation}\label{SkewSymetricPairing}
  \begin{aligned}
    & \lbrack (\bu,\bp),(\bv,\bq) \rbrack:=\langle \bu,\bq\rangle - \langle \bv,\bp\rangle\\
    & \bu,\bv\in \mbH(\Sigma),\;\bp,\bq\in \mbH(\Sigma)^{*}.
  \end{aligned}
\end{equation}
This form is obviously non-degenerate and can be used as a duality pairing
over the space of tuples of Dirichlet-Neumann pairs of traces. Indeed denote
\begin{equation*}
  \begin{aligned}
    & \mathscr{H}(\Sigma):=\mbH(\Sigma)\times\mbH(\Sigma)^{*}\\
    & \text{with norm:}\\
    & \Vert (\bv,\bq)\Vert_{\mT\times\mT^{-1}}^{2} :=
    \Vert \bv\Vert_{\mT}^{2}+\Vert \bq\Vert_{\mT^{-1}}^{2}
  \end{aligned}
\end{equation*}
then for any $\varphi\in \mathscr{H}(\Sigma)^{*}$,
there exists a unique $\ctru\in \mathscr{H}(\Sigma)$ such that 
$\lbrack \ctru,\ctrv \rbrack = \varphi(\ctrv)
\;\forall \ctrv\in \mathscr{H}(\Sigma)$. In other words,
the pairing \eqref{SkewSymetricPairing} puts $\mathscr{H}(\Sigma)$
in self-duality. We now introduce the subspace of so-called Cauchy
data that directly relates to the boundary value problem under study,
\begin{equation}\label{CauchyDataSpace}
  \mathscr{C}(\mA):=\{ (\mB(\bu),\bp)\;\vert\;
  (\bu,\bp)\in \mbH(\Omega)\times \mbH(\Sigma)^*,\;\mA\bu = \mB^{*}\bp\}
\end{equation}
It must be understood as the space of tuples of Dirichlet-Neumann trace pairs
stemming from solutions to the problems local to each subdomain.
If $\mA:\mbH(\Omega)\to \mbH(\Omega)^{*}$ is an isomorphism,
we can define the associated Neumann-to-Dirichlet operator
$\mrm{NtD}_{\mA}:=\mB\mA^{-1}\mB^{*}$ and then
$\mathscr{C}(\mA):=\{(\mrm{NtD}_{\mA}(\bp),\bp)\,\vert\,\bp\in\mbH(\Sigma)^{*}\}$
appears to be the graph of it. On the other hand $\mathscr{C}(\mA)$ is properly
defined even if $\mA$ fails to be invertible.

\begin{lem}\label{EnergyConservationIdentity}\quad\\
  Assume \eqref{UnifiedVariationalSetting}-\eqref{SuperBilinearForm}-\eqref{AssumptionAbsorption}-\eqref{LocalUniqueSolvability}.
  The application $(\bv,\bp)\to \bp-i\mT(\bv)$
  continuously and isomorphically maps $\mathscr{C}(\mA)$
  into $\mbH(\Sigma)^*$ and, for all  $(\bv,\bp)\in\mathscr{C}(\mA)$,
  satisfies the estimates
  \begin{equation*}
    \begin{aligned}
      \Vert \bv\Vert_{\mT}^{2} + \Vert \bp\Vert_{\mT^{-1}}^{2}
      & \leq \Vert \bp-i\mT\bv\Vert_{\mT^{-1}}^{2}\\
      \frac{1}{2}\Vert \bp-i\mT\bv\Vert_{\mT^{-1}}^{2}
      & \leq \Vert \bv\Vert_{\mT}^{2} + \Vert \bp\Vert_{\mT^{-1}}^{2}. 
    \end{aligned}
  \end{equation*}
\end{lem}
\noindent \textbf{Proof:}

It suffices to prove surjectivity and the estimates.
To prove surjectivity, pick an arbitrary $\bq\in\mbH(\Sigma)^{*}$
and define $\bu = (\mA-i\mB^*\mT\mB)^{-1}\mB^{*}\bq$. The pair
$(\bv,\bp) = (\mB(\bu), \bq+i\mT\mB(\bu))$ satisfies $\mA\bu = \mB^*\bp$
so that $(\bv,\bp)\in\mathscr{C}(\mA)$ and, by construction, we have
$\bp-i\mT\bv = \bq$.

To prove the estimates, pick an arbitrary pair $(\bv,\bp)\in\mathscr{C}(\mA)$.
According to \eqref{CauchyDataSpace} there exists $\bu\in \mbH(\Omega)$
such that $\mB(\bu) = \bv$ and $\mA(\bu) = \mB^*(\bp)$, hence
$\langle \bp,\overline{\bv}\rangle =
\langle \bp,\mB(\overline{\bu})\rangle = \langle \mB^*(\bp),\overline{\bu}\rangle
= \langle \mA(\bu),\overline{\bu}\rangle$. Taking account of
\eqref{GlobalAbsorptionProperty}, we deduce
$0\leq \Re e\{i\langle \bp,\overline{\bv}\rangle\}\leq
\Vert \bv\Vert_{\mT}^{2} + \Vert \bp\Vert_{\mT^{-1}}^{2}$ and conclude 
$0\leq \Vert \bp-i\mT\bv\Vert_{\mT^{-1}}^{2}  -(\Vert \bv\Vert_{\mT}^{2} + \Vert \bp\Vert_{\mT^{-1}}^{2})
\leq \Vert \bv\Vert_{\mT}^{2} + \Vert \bp\Vert_{\mT^{-1}}^{2}$. 
\hfill $\Box$

\quad\\
In the previous lemma, the space of Cauchy data has been proven boundedly
isomorphic to a Hilbert space and, as such, is closed. 

\begin{cor}\quad\\
  Assume \eqref{UnifiedVariationalSetting}-\eqref{SuperBilinearForm}-\eqref{AssumptionAbsorption}-\eqref{LocalUniqueSolvability}.
  The subspace $\mathscr{C}(\mA)$ is closed in $\mathscr{H}(\Sigma)$.
\end{cor}

\noindent 
The space of Cauchy data can be complemented in various ways. The next proposition
exhibits one possibility.

\begin{prop}\label{DirectSumDirichletNeumannPairs}\quad\\
  Assume \eqref{UnifiedVariationalSetting}-\eqref{SuperBilinearForm}-\eqref{AssumptionAbsorption}-\eqref{LocalUniqueSolvability}.
  Define  $\mathscr{G}(i\mT):=\{ (\bv,i\mT(\bv)), \bv\in \mbH(\Sigma)\}$. Then
  \begin{equation*}
    \mathscr{H}(\Sigma) = \mathscr{C}(\mA)\oplus \mathscr{G}(i\mT).
  \end{equation*}
\end{prop}
\noindent\textbf{Proof:}

First of all, assume that $(\bu,\bp)\in \mathscr{C}(\mA)\cap \mathscr{G}(i\mT)$.
This that there exists $\bv\in \mbH(\Omega)$ such that $\mA\bv = \mB^{*}\bp$ and
$\mB\bv = \bu$, and that $\bp = i\mT\bu$. Combining these equations yields
$(\mA-i\mB^*\mT\mB)\bv = 0$ hence $\bv=0$ according to Lemma \ref{MultiDomainAbsorptionProp},
and finally $(\bu,\bp) = 0$. We have proved that $\mathscr{C}(\mA)\cap \mathscr{G}(i\mT) = \{0\}$.

Now take an arbitrary $(\bu,\bp)\in \mbH(\Sigma)\times\mbH(\Sigma)^{*}$.
Since $\mB:\mbH(\Omega)\to\mbH(\Sigma)$ is surjective, there exists $\bw\in \mbH(\Omega)$
such that $\mB(\bw) = \bu$. Define $\bv\in \mbH(\Omega)$ by 
$\bv = (\mA-i\mB^*\mT\mB)^{-1}(\mA\bw-\mB^{*}\bp)$ which is valid a definition
since $\mA-i\mB^*\mT\mB:\mbH(\Omega)\to \mbH(\Omega)^{*}$ is an isomorphism
according to Lemma \ref{MultiDomainAbsorptionProp}. We have in particular
$\mA(\bw-\bv) = \mB^{*}(\bp-i\mT\mB\bv)$. Set
\begin{equation}\label{DecompositionGraph}
  \begin{aligned}
    & \bu_1 = \mB(\bv),\quad \bp_1 = i\mT\bu_1  = i\mT\mB(\bv),\\
    & \bu_2 = \mB(\bw-\bv) = \bu - \bu_1,\quad \bp_2 = \bp-i\mT\mB\bv = \bp-\bp_1.
  \end{aligned}
\end{equation}
By construction we have $(\bu_1,\bp_1)\in\mathscr{G}(i\mT)$. Moreover
$\mB(\bw-\bv) = \bu_2$ and $\mA(\bw-\bv) = \mB^{*}\bp_2$ so that
$(\bu_2,\bp_2)\in\mathscr{C}(\mA)$. Finally, the second line in
\eqref{DecompositionGraph} indicates that
$(\bu,\bp) = (\bu_1,\bp_1) +  (\bu_2,\bp_2)$ which thus proves 
$(\bu,\bp)\in \mathscr{C}(\mA)+\mathscr{G}(i\mT)$. We have just established that 
$\mathscr{C}(\mA)+\mathscr{G}(i\mT) = \mbH(\Sigma)\oplus\mbH(\Sigma)^{*}$ which ends
the proof. \hfill $\Box$

\quad\\
The space $\mathscr{G}(i\mT)$ is simply the graph of the (bounded) operator
$i\mT:\mbH(\Sigma)\to \mbH(\Sigma)^{*}$. In the present analysis, it plays a
secondary role and shall be used only to prove results about $\mathscr{C}(\mA)$. 
We have the following immediate result.

\begin{lem}\label{SelfPolarityGraphiT}\quad\\
  Define $\mathscr{G}(i\mT)^{\sharp} := \{\bu\in \mathscr{H}(\Sigma),\;
  \lbrack\bu,\bv\rbrack = 0\; \forall\bv\in \mathscr{G}(i\mT)\}$. Then
  $\mathscr{G}(i\mT)^{\sharp} = \mathscr{G}(i\mT)$.
\end{lem}

\noindent 
The proof is definitely straightforward. This result means that $\mathscr{G}(i\mT)$
is its own polar set under the pairing $\lbr\cdot,\cdot\rbr$. As we see now, the space
$\mathscr{C}(\mA)$ fulfills a similar property.

\begin{prop}\label{CauchyTracesSelfPolarity}\quad\\
  Assume \eqref{UnifiedVariationalSetting}-\eqref{SuperBilinearForm}-\eqref{AssumptionAbsorption}-\eqref{LocalUniqueSolvability}.
  Define $\mathscr{C}(\mA)^{\sharp} := \{\bu\in \mathscr{H}(\Sigma),\;
  \lbrack\bu,\bv\rbrack = 0\; \forall\bv\in \mathscr{C}(\mA)\}$. Then
  \begin{equation*}
    \mathscr{C}(\mA)^{\sharp} = \mathscr{C}(\mA^{*}).
  \end{equation*}
\end{prop}
\noindent\textbf{Proof:}

First of all we have $\mathscr{C}(\mA^{*})\subset \mathscr{C}(\mA)^{\sharp}$.
Indeed take any $(\bu,\bp)\in\mathscr{C}(\mA)$. By definition, there exists
$\bw\in \mbH(\Omega)$ such that $\mB(\bw) = \bu$ and $\mA\bw = \mB^*\bp$.
Then for any $(\bu',\bp')\in\mathscr{C}(\mA^*)$, since
$\mB(\bw') = \bu'$ and $\mA^*\bw' = \mB^*\bp'$ for some
$\bw'\in\mbH(\Omega)$, we have
\begin{equation*}
  \begin{aligned}
    \lbr (\bu,\bp), (\bu',\bp')\rbr
    & = \langle \bu,\bp'\rangle - \langle \bu',\bp\rangle =
    \langle \mB(\bw),\bp'\rangle - \langle \mB(\bw'),\bp\rangle \\
    & = \langle \bw,\mB^*(\bp')\rangle - \langle \bw',\mB^*(\bp)\rangle\\
    & = \langle \bw,\mA^*(\bw')\rangle - \langle \bw',\mA(\bw)\rangle = 0.
  \end{aligned}
\end{equation*}
Hence, to finish the proof, we need to show that
$\mathscr{C}(\mA)^{\sharp}\subset \mathscr{C}(\mA^{*})$. For that,
pick an arbitrary $\ctru = (\bu,\bp)\in \mathscr{C}(\mA)^{\sharp}$.
The hypothesis of Section \ref{SingleDomainVarForm} hold for
$\mA_{\Omega\times\Gamma}^{*}$ instead of $\mA_{\Omega\times\Gamma}$,
hence we can apply Proposition \ref{DirectSumDirichletNeumannPairs}
to $\mA^*$. This yields a decomposition
$\ctru = \ctru_1 + \ctru_2$ for some $\ctru_1\in
\mathscr{C}(\mA^*)$ and some $\ctru_2\in\mathscr{G}(i\mT)$. We have
to prove that $\ctru_2 = 0$. By assumption we have
\begin{equation*}
  0 = \lbr \ctru,\ctrv\rbr =   \lbr \ctru_1,\ctrv\rbr +   \lbr \ctru_2,\ctrv\rbr =   \lbr \ctru_2,\ctrv\rbr
  \quad\forall \ctrv\in \mathscr{C}(\mA),
\end{equation*}
since $\mathscr{C}(\mA)\subset \mathscr{C}(\mA^*)^{\sharp}$.
Next Lemma \ref{SelfPolarityGraphiT} implies that $0 = \lbr\ctru_2,\ctrv\rbr =
\lbr\ctru_2,\ctrv+\ctrv'\rbr$ for all $\ctrv\in \mathscr{C}(\mA)$ and all
$\ctrv'\in \mathscr{G}(i\mT)$. Since $\mathscr{C}(\mA)\oplus\mathscr{G}(i\mT) = \mathscr{H}(\Sigma)$
according to Proposition \ref{DirectSumDirichletNeumannPairs}, we conclude that
$0 = \lbr\ctru_2,\ctrw\rbr\;\forall \ctrw\in  \mathscr{H}(\Sigma)$
hence finally $\ctru_2=0$. This shows that $\ctru = \ctru_1\in\mathscr{C}(\mA^*)$.
We have just established that $\mathscr{C}(\mA)^{\sharp}\subset \mathscr{C}(\mA^*)$.
\hfill $\Box$

\quad\\
We point that, because $\mathscr{C}(\mA)$ is closed, the previous result also
implies that $\mathscr{C}(\mA) = \mathscr{C}(\mA^*)^{\sharp}$. 
Self-polarity appears to be a property of the following subspace (see Proposition
\ref{NonlocalTrCond}) that is pivotal in characterizing transmission conditions
\begin{equation*}
  \mathscr{X}(\Sigma):=\mbX(\Sigma)\times\mbX(\Sigma)^{\circ}.
\end{equation*}
Indeed we have $\mathscr{X}(\Sigma) =\mathscr{X}(\Sigma)^{\sharp} := \{\ctru\in
\mathscr{H}(\Sigma),\;\lbr\ctru,\ctrv\rbr = 0\;\forall\ctrv\in \mathscr{X}(\Sigma)\}$ by
the very definition of $\mathscr{X}(\Sigma)$, as $\mbX(\Sigma)^{\circ\circ} = \mbX(\Sigma)$
since $\mbX(\Sigma)$ is a closed subspace of $\mbH(\Sigma)$ (see e.g.
\cite[Thm.4.7]{zbMATH01022519} or \cite[Prop.1.9]{MR2759829}). The next result establishes an
important connection between the two spaces $\mathscr{C}(\mA),\mathscr{X}(\Sigma)$ and
our primary boundary value problem \eqref{PrimaryFormulation}.

\begin{prop}\label{CorrespondanceKernelIntersection}\quad\\
  Assume \eqref{UnifiedVariationalSetting}-\eqref{SuperBilinearForm}-\eqref{AssumptionAbsorption}-\eqref{LocalUniqueSolvability}.
  The operator $\bu \mapsto (\mB\mR(\bu), (\mB^{\dagger})^*\mA\mR(\bu))$
  continuously and isomorphically maps $\mrm{ker}(\mA_{\Omega\times\Gamma})$
  onto $\mathscr{C}(\mA)\cap\mathscr{X}(\Sigma)$. As a consequence
  \begin{equation*}
    \mrm{dim}(\mrm{ker}(\mA_{\Omega\times\Gamma})\,) = \mrm{dim}(\mathscr{C}(\mA)\cap \mathscr{X}(\Sigma)).
  \end{equation*}
\end{prop}
\noindent \textbf{Proof:}

Let $\bu\in \mH(\Omega\times\Gamma)$ satisfy $\mA_{\Omega\times\Gamma}(\bu) = 0$.
In particular $\mR(\bu)\in\mbX(\Omega)$ and $\mA\mR(\bu)\in \mbX(\Omega)^{\circ}$,
see \eqref{DefOperatorRestriction} and \eqref{DomainDecomposedPb1}.
According to \textit{iv)} of Lemma \ref{ExtendedPolarity}, there exists
$\bp\in \mbX(\Sigma)^{\circ}$ such that $\mA\mR(\bu) = \mB^{*}\bp$ and it is
unique since $\mB^*:\mbH(\Sigma)^{*}\to \mbH(\Omega)^{*}$ is injective.
We have
\begin{equation*}
  (\mB^{\dagger})^{*}\mA\mR(\bu) = (\mB^{\dagger})^{*}\mB^*\bp = (\mB\mB^{\dagger})^{*}\bp = \bp.
\end{equation*}
Setting $\bv := \mB\cdot\mR(\bu)$, by construction $(\bv,\bp)\in\mathscr{C}(\mA)$.
We also have $\bv\in \mbX(\Sigma)$ since $\mR(\bu)\in\mbX(\Omega)$,
so that $(\bv,\bp)\in \mbX(\Sigma)\times \mbX(\Sigma)^{\circ} = \mathscr{X}(\Sigma)$.
In addition, the formula $(\bv,\bp) = (\mB\mR\bu, (\mB^{\dagger})^*\mA\mR\bu)$ establishes
the continuous dependency of $(\bv,\bp)$ on $\bu$.

Reciprocally, consider an arbitrary pair $(\bv,\bp)\in \mathscr{C}(\mA)\cap \mathscr{X}(\Sigma)$.
Since $(\bv,\bp)\in \mathscr{C}(\mA)$, there exists $\bw\in \mbH(\Omega)$ such that
$\mA\bw = \mB^{*}\bp$ and $\mB(\bw) = \bv$, and such a $\bw$ is unique since
$\mrm{ker}(\mA)\cap \mrm{ker}(\mB) = \{0\}$, according to Lemma \ref{MultiDomainAbsorptionProp}.
As $\bv\in \mbX(\Sigma)$, we have $\bw\in \mbX(\Omega) = \mB^{-1}(\mbX(\Sigma)) $
according to \textit{iii)} of Lemma \ref{ExtendedPolarity}, so there exists $\bu\in
\mH(\Omega\times\Gamma)$ such that $\mR(\bu) = \bw$ and such a $\bu$ is unique due
to the injectivity of $\mR:\mH(\Omega\times\Gamma)\to \mbH(\Omega)$. This leads to 
$\mA\mR(\bu) = \mB^{*}\bp$ and $\bp\in \mbX(\Sigma)^{\circ}\Rightarrow \mB^{*}\bp\in
\mbX(\Omega)^{\circ} = \mrm{ker}(\mR^*)$. Since $\mbX(\Omega) = \mR(\mH(\Omega\times\Gamma))$, we
conclude that $0 = \mR^{*}\mA\mR(\bu) = \mA_{\Omega\times\Gamma}(\bu)$.  \hfill $\Box$

\begin{lem}\label{SurjectivityPreparatoryLemma}\quad\\
  Assume \eqref{UnifiedVariationalSetting}-\eqref{SuperBilinearForm}-\eqref{AssumptionAbsorption}-\eqref{LocalUniqueSolvability}.
  The operator $(\bu,\bp)\mapsto \mR^*(\mB^*\bp - \mA\mB^{\dagger}\bu)$
  continuously maps $(\mathscr{C}(\mA^{*})\cap\mathscr{X}(\Sigma))^{\sharp}$
  into $\mrm{range}(\mA_{\Omega\times\Gamma})$.
\end{lem}
\noindent \textbf{Proof:}

Take an arbitrary $(\bu,\bp)\in(\mathscr{C}(\mA^{*})\cap\mathscr{X}(\Sigma))^{\sharp}$
and set $\bff = \mR^*(\mB^*\bp - \mA\mB^{\dagger}\bu)$. Applying Proposition
\ref{CorrespondanceKernelIntersection} to $\mA_{\Omega\times\Gamma}^{*}$
instead of $\mA_{\Omega\times\Gamma}$ shows that $\varphi\in \mrm{ker}(\mA_{\Omega\times\Gamma}^{*})
\Rightarrow (\bv,\bq) = (\mB\mR(\varphi), (\mB^{\dagger})^*\mA^*\mR(\varphi)) \in
\mathscr{C}(\mA^*)\cap \mathscr{X}(\Sigma)$. Hence
$\langle \bff,\varphi\rangle
= \langle \mR^*(\mB^*\bp - \mA\mB^{\dagger}\bu),\varphi\rangle
= \langle \bp, \mB\mR\varphi\rangle - \langle \bu,(\mB^\dagger)^*\mA^*\mR\varphi\rangle
= \lbrack (\bv,\bq), (\bu,\bp) \rbrack = 0$.
This proves $\bff\in \mrm{ker}(\mA_{\Omega\times\Gamma}^*)^{\circ}
= \mrm{range}(\mA_{\Omega\times\Gamma})$ according to \eqref{FredholmAlternative}.
\hfill $\Box$

\begin{prop}\label{ClosedImage}\quad\\
  Assume \eqref{UnifiedVariationalSetting}-\eqref{SuperBilinearForm}-\eqref{AssumptionAbsorption}-\eqref{LocalUniqueSolvability}. Then $\mathscr{C}(\mA)+\mathscr{X}(\Sigma) = (\mathscr{C}(\mA^*)\cap\mathscr{X}(\Sigma))^{\sharp}$.
  In particular the subspace $\mathscr{C}(\mA)+\mathscr{X}(\Sigma)$ is closed in $\mathscr{H}(\Sigma)$.
\end{prop}
\noindent \textbf{Proof:}

Clearly we have $\mathscr{C}(\mA)+\mathscr{X}(\Sigma)\subset
(\mathscr{C}(\mA^*)\cap\mathscr{X}(\Sigma))^{\sharp}$, so we only need to
establish that $(\mathscr{C}(\mA^*)\cap\mathscr{X}(\Sigma))^{\sharp}\subset
\mathscr{C}(\mA)+\mathscr{X}(\Sigma)$. Pick any pair $(\bp_{\dir},\bp_{\neu})
\in (\mathscr{C}(\mA^*)\cap\mathscr{X}(\Sigma))^{\sharp}$. According to
Lemma \ref{SurjectivityPreparatoryLemma} we have
$\mR^*(\mB^*\bp_{\neu} - \mA\mB^{\dagger}\bp_{\dir})\in\mrm{range}(\mA_{\Omega\times\Gamma})$.
Applying the definition of $\mA$ given by
\eqref{DefinitionOperatorMultiDomain}, there exists
$\varphi \in\mbX(\Omega)$ satisfying $\langle\mA\varphi,\bw\rangle =
\langle \mB^*\bp_{\neu}-\mA\mB^{\dagger}\bp_{\dir},\bw\rangle$ for all
$\forall \bw\in \mbX(\Omega)$.

Set $\phi  = \varphi+\mB^{\dagger}(\bp_{\dir})$
and $\bu_{\dir} = \mB(\phi) = \mB(\varphi)+\bp_{\dir}$. By construction,
$\langle \mA(\phi),\bw\rangle = \langle\bp_{\neu},\mB(\bw)\rangle  = 0
\;\forall \bw\in \mrm{ker}(\mB)\subset \mbX(\Omega)$, which rewrites
$\mA(\phi)\in \mrm{ker}(\mB)^{\circ}$. Applying \textit{i)}
of Lemma \ref{ExtendedPolarity} we have $\mA\phi = \mB^*\bu_{\neu}$ for
some $\bu_{\neu}\in \mbH(\Sigma)^{*}$. This implies in particular
$\bu_{\neu} = (\mB\mB^{\dagger})^*\bu_{\neu} = (\mB^{\dagger})^*\mB^*\bu_{\neu} =
(\mB^{\dagger})^*\mA\phi$.

We have $\mA\phi = \mB^*\bu_{\neu}$ and $\mB\phi = \bu_{\dir}$ hence
$(\bu_{\dir},\bu_{\neu})\in \mathscr{C}(\mA)$.
On the other hand $\bp_{\dir}-\bu_{\dir} = -\mB\varphi\in \mbX(\Sigma)$
since $\varphi\in\mbX(\Omega)$ and, for any $\bw\in \mbX(\Sigma)$ we have
$\mB^{\dagger}(\bw)\in\mbX(\Omega)$ hence 
$\langle \bp_{\neu}-\bu_{\neu},\bw\rangle  =
\langle \mA\phi,\mB^{\dagger}\bw\rangle - \langle \mA\phi,\mB^{\dagger}\bw\rangle = 0$,
which implies $\bp_{\neu}-\bu_{\neu}\in\mbX(\Sigma)^{\circ}$.
Finally $(\bu_{\dir},\bu_{\neu})\in\mathscr{C}(\mA)$
and $(\bp_{\dir},\bp_{\neu}) - (\bu_{\dir},\bu_{\neu})\in \mathscr{X}(\Sigma)$
imply that $(\bp_{\dir},\bp_{\neu})\in \mathscr{C}(\mA)+ \mathscr{X}(\Sigma)$.
\hfill $\Box$

\begin{cor}\label{CokernelCorrespondance}\quad\\
  Assume \eqref{UnifiedVariationalSetting}-\eqref{SuperBilinearForm}-\eqref{AssumptionAbsorption}-\eqref{LocalUniqueSolvability}.
  Then
  \begin{equation*}
    \mrm{codim}(\mathscr{C}(\mA)+\mathscr{X}(\Sigma)\,) =  
    \mrm{codim}(\mrm{range}(\mA_{\Omega\times\Gamma})\,).
  \end{equation*}
\end{cor}
\noindent \textbf{Proof:}

We have $(\mathscr{C}(\mA)+\mathscr{X}(\Sigma))^{\sharp} =
\mathscr{C}(\mA)^{\sharp}\cap \mathscr{X}(\Sigma)^{\sharp}$ see
e.g. \cite[Prop.2.14]{MR2759829}. According to Proposition
\ref{CauchyTracesSelfPolarity} applied to $\mA^*$, and since
$\mathscr{X}(\Sigma)^{\sharp} = \mathscr{X}(\Sigma)$ by construction,
we conclude that $(\mathscr{C}(\mA)+\mathscr{X}(\Sigma))^{\sharp} =
\mathscr{C}(\mA^*)\cap\mathscr{X}(\Sigma)$. As the bilinear pairing
$\lbr\cdot,\cdot\rbr$ is non-degenerate and
$\mathscr{C}(\mA)+\mathscr{X}(\Sigma)$ is closed according to
Proposition \ref{ClosedImage}, we conclude
$\mrm{codim}(\mathscr{C}(\mA)+\mathscr{X}(\Sigma)) =
\mrm{dim}((\mathscr{C}(\mA)+\mathscr{X}(\Sigma))^{\sharp}) =
\mrm{dim}(\mathscr{C}(\mA^*) \cap\mathscr{X}(\Sigma))$.
There only remains to apply Proposition \ref{CorrespondanceKernelIntersection}
to $\mA_{\Omega\times\Gamma}^{*}$ combined with \eqref{FredholmAlternative}.
\hfill $\Box$

\section{Scattering operator}

Proposition \ref{CorrespondanceKernelIntersection} and \ref{ClosedImage} and
Corollary \ref{CokernelCorrespondance} above show that
the kernel and the range of $\mA_{\Omega\times\Gamma}$ are closely related to 
the pair of subspaces $\mathscr{C}(\mA),\mathscr{X}(\Sigma)$. This can be
exploited to study other formulations of the same boundary value problem.

\begin{prop}\quad\\
  Assume \eqref{UnifiedVariationalSetting}-\eqref{SuperBilinearForm}-\eqref{AssumptionAbsorption}-\eqref{LocalUniqueSolvability}.
  If $\bu\in \mbX(\Omega)$ satisfies \eqref{DomainDecomposedPb1},
  then there exists a unique $\bp\in \mbH(\Sigma)^{*}$ such that
  the pair $(\bu,\bp)$ satisfies
  \begin{equation}\label{DomainDecomposedPb2}
    \begin{aligned}
      & \bu\in\mbH(\Omega),\;\bp\in \mbH(\Sigma)^*,\\
      & \mA\bu - \mB^*\bp = \bell,\\
      & -\bp + i\mT\mB\bu = \Pi(\bp+i\mT\mB\bu).
    \end{aligned}
  \end{equation}
  Reciprocally if the pair $(\bu,\bp)\in\mbH(\Omega)\times\mbH(\Sigma)^{*}$ satisfies
  \eqref{DomainDecomposedPb2}, then $\bu$ satisfies \eqref{DomainDecomposedPb1}.
\end{prop}
\noindent \textbf{Proof:}

Assume first that $\bu\in \mbX(\Omega)$ satisfies \eqref{DomainDecomposedPb1}.
This formulation rewrites equivalently as $\mA\bu - \bell\in\mbX(\Omega)^{\circ}$.
Since $\mbX(\Omega)^{\circ} = \mB^{*}(\mbX(\Sigma)^{\circ})$ according to \textit{iv)}
Lemma \ref{ExtendedPolarity}, and as $\mB^{*}:\mbH(\Sigma)^*\to \mbH(\Omega)^*$ is
injective ($\mB$ is surjective), there exists a unique $\bp\in \mbX(\Sigma)^\circ$
such that $\mA\bu - \bell = \mB^*\bp$. On the other hand, $\bu \in\mbX(\Omega)\Rightarrow \mB(\bu)
\in\mbX(\Sigma)$ according to \textit{iii)} of Lemma \ref{ExtendedPolarity}. Finally
applying Proposition \ref{NonlocalTrCond}, we obtain $-\bp + i\mT\mB\bu = \Pi(\bp+i\mT\mB\bu)$.

Reciprocally, assume that \eqref{DomainDecomposedPb2} holds. Then, according to
Proposition \ref{NonlocalTrCond}, we have $\bp\in \mbX(\Sigma)^{\circ}$ and
$\mB(\bu)\in \mbX(\Sigma)$. Moreover we have
$\mB(\bu)\in \mbX(\Sigma)\Rightarrow \bu\in\mbX(\Omega)$
according to \textit{iii)} of Lemma \ref{ExtendedPolarity}.
Since $\bp\in \mbX(\Sigma)^{\circ}$, we have $\mB^{*}\bp\in\mbX(\Omega)^{\circ}$
so that, for any $\bv\in\mbX(\Omega)$ we have $0 = \langle \mB^{*}\bp,\bv\rangle =
\langle \mA\bu-\bell,\bv\rangle$. To sum up, we have proved that $\bu\in\mbX(\Omega)$
and $\langle \mA\bu,\bv\rangle = \langle \bell,\bv\rangle\;\forall\bv\in\mbX(\Omega)$.
\hfill $\Box$

\quad\\
In a domain decomposition context, a substructuring strategy
applied to Problem \eqref{PrimaryFormulation} naturally leads
to eliminating the volume unknowns in \eqref{DomainDecomposedPb2}.
This is performed by means of a scattering map that takes ingoing
traces as input and returns outgoing traces as output.

\begin{prop}\label{DefScatOp2}\quad\\
  Assume \eqref{UnifiedVariationalSetting}-\eqref{SuperBilinearForm}-\eqref{AssumptionAbsorption}-\eqref{LocalUniqueSolvability}.
  There exists a unique bounded linear map $\mS:\mbH(\Sigma)^{*}\to \mbH(\Sigma)^{*}$,
  later referred to as scattering operator,  satisfying
  \begin{equation}\label{DefScatteringOperator}
    \bp+i\mT\bv = \mS(\bp-i\mT\bv)\quad \forall (\bv,\bp)\in\mathscr{C}(\mA).
  \end{equation}
  It is also given by the formula $\mS=\Id +2i\mT\mB(\mA-i\mB^*\mT\mB)^{-1}\mB^*$.
  It is $\mT^{-1}$-contractive and, for any $\bq\in \mbH(\Sigma)^{*}$, satisfies  
  \begin{equation*}
    \begin{aligned}
    & \Vert \mS(\bq)\Vert_{\mT^{-1}}^{2} +4\vert\Im m\{\langle
    \mA(\bu),\overline{\bu}\rangle\}\vert = \Vert \bq\Vert_{\mT^{-1}}^{2}\\
    & \text{where}\;\; \bu  = (\mA-i\mB^*\mT\mB)^{-1}\mB^*\bq.
    \end{aligned}
  \end{equation*}
\end{prop}
\noindent \textbf{Proof:}

We follow the proof pattern presented e.g. in \cite[Lem.5.2]{claeys2021nonself}. 
First of all, Identity \eqref{DefScatteringOperator} clearly and unambiguously defines
the operator $\mS$ as a linear map according to Lemma \ref{EnergyConservationIdentity}.
Next, pick an arbitrary $\bq\in \mbH(\Sigma)^{*}$ and set $\bu = (\mA-i\mB^{*}\mT\mB)^{-1}\mB^{*}\bq$
and $\bp = \bq+i\mT\mB(\bu)$. We have $\mA\bu-\mB^*\bp = 0$ and $\bq = \bp - i\mT\mB(\bu)$
and $\mS(\bq) = \bp + i\mT\mB(\bu) = \bq+2i\mT\mB(\bu)$, which leads to
$\mS(\bq)=(\Id +2i\mT\mB(\mA-i\mB^*\mT\mB)^{-1}\mB^*)\bq$.
Finally developing the squared norm, and taking account of
\eqref{GlobalAbsorptionProperty}, we have
\begin{equation*}
  \begin{aligned}
    \Vert \mS(\bq)\Vert_{\mT^{-1}}^{2}
    & = \Vert \bp + i\mT\mB(\bu)\Vert_{\mT^{-1}}^{2}\\
    & = \Vert \bp - i\mT\mB(\bu)\Vert_{\mT^{-1}}^{2} + 4 \Im m\{\langle\bq,\mB(\overline{\bu})\rangle\}
    +4\Vert\mB(\bu)\Vert_{\mT}^{2}\\
    & = \Vert \bq\Vert_{\mT^{-1}}^{2} + 4 \Im m\{\langle\mB^*(\bq),\overline{\bu}\rangle\}
    +4\Vert\mB(\bu)\Vert_{\mT}^{2}\\
    & = \Vert \bq\Vert_{\mT^{-1}}^{2} + 4 \Im m\{\langle\mA(\bu),\overline{\bu}\rangle\}
    - 4\Im m\{i\langle \mB^*\mT\mB(\bu),\overline{\bu}\rangle  \}+4\Vert\mB(\bu)\Vert_{\mT}^{2}\\
    & =  \Vert \bq\Vert_{\mT^{-1}}^{2} - 4 \vert \Im m\{\langle\mA(\bu),\overline{\bu}\rangle\}\vert
  \end{aligned}
\end{equation*}
\hfill $\Box$

\quad\\
The space of Cauchy data was used to characterize the scattering operator.
Reciprocally, the scattering operator provides a characterization of the space
of Cauchy data. The following result should be compared with
\eqref{CaracSingleTraceSpace}.

\begin{lem}\label{CaracScatteringOperator}\quad\\
  Assume \eqref{UnifiedVariationalSetting}-\eqref{SuperBilinearForm}-\eqref{AssumptionAbsorption}-\eqref{LocalUniqueSolvability}.
  For any $(\bv,\bp)\in\mathscr{H}(\Sigma)$ we have:
  \begin{equation*}
    (\bv,\bp)\in\mathscr{C}(\mA)\iff\bp+i\mT\bv = \mS(\bp-i\mT\bv).
  \end{equation*}
\end{lem}
\noindent \textbf{Proof:}

From the very definition of the scattering operator in Proposition \ref{DefScatOp2},
it is clear that $(\bv,\bp)\in\mathscr{C}(\mA)\Rightarrow\bp+i\mT\bv =
\mS(\bp-i\mT\bv)$. Reciprocally pick arbitrarily some $(\bv,\bp)\in\mathscr{H}(\Sigma)$
such that $\bp+i\mT\bv = \mS(\bp-i\mT\bv)$. We know from Proposition
\ref{DirectSumDirichletNeumannPairs} that there exists
$\bv'\in\mbH(\Sigma)$ such that $(\bv-\bv',\bp-i\mT\bv')\in\mathscr{C}(\mA)$
so applying Proposition \ref{DefScatOp2} we obtain
\begin{equation*}
  \begin{aligned}
    & (\bp-i\mT\bv') + i\mT(\bv-\bv')
    = \mS(\, (\bp-i\mT\bv') - i\mT(\bv-\bv')\,)\\
    & \iff\;\; \bp + i\mT\bv - 2i\mT\bv'
    = \mS(\bp-i\mT\bv) \\
    & \iff\;\; 2i\mT\bv' = 0\;\; \Longrightarrow \;\; \bv' = 0.
  \end{aligned}
\end{equation*}
\hfill $\Box$

\quad\\
The scattering operator has a subdomain-wise block diagonal structure.
This is clearly visible from the formula $\mS=\Id +2i\mT\mB(\mA-i\mB^*\mT\mB)^{-1}\mB^*$
where each term in the right hand side is block diagonal. This yields
\begin{equation*}
  \begin{aligned}
    & \mS = \mrm{diag}(\mS_{\Gamma},\mS_{\Omega_1},\dots,\mS_{\Omega_{\mJ}})\\
    & \text{where}\;\; \mS_{\Omega_j} \!= \Id +2i\mT_{\Omega_j}\mB_{\Omega_j}
    (\mA_{\Omega_j}-i\mB_{\Omega_j}^*\mT_{\Omega_j}\mB_{\Omega_j})^{-1}\mB_{\Omega_j}^*\\
    & \textcolor{white}{where}\;\; \mS_{\Gamma} = \Id +2i\mT_{\Gamma}\mB_{\Gamma}(
    \mA_{\Gamma}-i\mB_{\Gamma}^*\mT_{\Gamma}\mB_{\Gamma})^{-1}\mB^*_{\Gamma}
  \end{aligned}
\end{equation*}
Let us discuss the particular form that takes the boundary scattering
operator $\mS_{\Gamma}$ for Dirichlet, Neumann and Robin conditions.
Recall that $\mB_{\Gamma}:\mH^{\textsc{b}}(\Gamma) := \mH^{1/2}(\Gamma)
\times\mH^{-1/2}(\Gamma)\to \mH^{1/2}(\Gamma)$ is defined by
$\mB_{\Gamma}(\alpha,p) = \alpha$ hence $\mB_{\Gamma}^*(p) = (p,0)$.

\begin{example}[\textbf{Dirichlet condition}]
  Taking the same notations as in Example \ref{DirichletBC1} and \ref{DirichletBC2},
  since $\mB_{\Gamma}^*p = (p,0)$ for all $p\in \mH^{-1/2}(\Gamma)$, we conclude that 
  $\mB_{\Gamma}(\mA_{\Gamma}-i\mB_{\Gamma}^*\mT_{\Gamma}\mB_{\Gamma})^{-1}\mB_{\Gamma}^*= 0$
  and finally
  \begin{equation*}
    \mS_{\Gamma} = +\Id.
  \end{equation*}
\end{example}

\begin{example}[\textbf{Neumann condition}]
  Taking the same notations as in Example \ref{NeumannBC1} and \ref{NeumannBC2},  
  in this situation we have $\mB_{\Gamma}(\mA_{\Gamma}-i\mB_{\Gamma}^*\mT_{\Gamma}
  \mB_{\Gamma})^{-1}\mB_{\Gamma}^*= i\mT_{\Gamma}^{-1}$. This yields the expression
  \begin{equation*}
    \mS_{\Gamma} = -\Id.
  \end{equation*}
\end{example}

\begin{example}[\textbf{Robin condition}]
  Taking the same notations as in Example \ref{RobinBC1} and \ref{RobinBC2},  
  in this situation we have $\mB_{\Gamma}(\mA_{\Gamma}-i\mB_{\Gamma}^*\mT_{\Gamma}
  \mB_{\Gamma})^{-1}\mB_{\Gamma}^*= i(\Lambda+\mT_{\Gamma})^{-1}$ which yields
  \begin{equation*}
    \mS_{\Gamma} = (\Lambda-\mT_{\Gamma})(\Lambda+\mT_{\Gamma})^{-1}.
  \end{equation*}  
\end{example}
  
\begin{example}[\textbf{Mixed conditions}]
Taking the same notations as in Example \ref{MixedBC1} and \ref{MixedBC2},
  in this situation we have   $\mB_{\Gamma}(\mA_{\Gamma}-i\mB_{\Gamma}^*\mT_{\Gamma}
  \mB_{\Gamma})^{-1}\mB_{\Gamma}^*= i\mT_{\Gamma}^{-1}(\Id-\Theta)$ which yields
  \begin{equation*}
    \mS_{\Gamma} = \Theta - (\Id-\Theta).
  \end{equation*}
  With this expression, the operator $\mS_{\Gamma}$ appears to be a
  $\mT_{\Gamma}^{-1}$-orthogonal symmetry with respect to the space
  $\widetilde{\mH}^{-1/2}(\Gamma_{\dir})$ and, in particular, a
  $\mT_{\Gamma}^{-1}$-isometry.
\end{example}

\section{Skeleton formulation}

Now we shall use the scattering operator of the previous section
to transform further the boundary value problem \eqref{DomainDecomposedPb2}.
Once volume unknowns have been eliminated, this reduces to an equation
involving only traces on the skeleton of the subdomain partition.

\begin{prop}\label{EquivalenceFormulation3}\quad\\
  Assume \eqref{UnifiedVariationalSetting}-\eqref{SuperBilinearForm}-\eqref{AssumptionAbsorption}-\eqref{LocalUniqueSolvability}.
  Define $\bff\in\mbH(\Sigma)^{*}$ by $\bff = -2i\Pi\mT\mB(\mA-i\mB^*\mT\mB)^{-1}\bell$.
  If $(\bu,\bp)\in \mbH(\Omega)\times \mbH(\Sigma)^{*}$ solves
  \eqref{DomainDecomposedPb2}, then $\bq = \bp-i\mT\mB(\bu)$
  satisfies the skeleton problem
  \begin{equation}\label{SkeletonEquation}
    \begin{aligned}
      & \bq\in \mbH(\Sigma)^*\;\text{and}\\
      & (\Id + \Pi\mS)\bq = \bff.
    \end{aligned}
  \end{equation}
  Reciprocally if $\bq$ satisfies the above equation then the pair
  $(\bu,\bp)\in \mbH(\Omega)\times \mbH(\Sigma)^{*}$, given by
  $\bu = (\mA-i\mB^*\mT\mB)^{-1}(\mB^*\bq + \bell)$
  and $\bp = \bq + i\mT\mB(\bu)$, solves \eqref{DomainDecomposedPb2}.
\end{prop}
\noindent \textbf{Proof:}

If $(\bu,\bp)\in \mbH(\Omega)\times \mbH(\Sigma)^{*}$ solves
\eqref{DomainDecomposedPb2} and $\bq = \bp-i\mT\mB(\bu)$, then
$(\mA-i\mB^*\mT\mB)\bu = \mB^{*}(\bp-i\mT\mB\bu)+\bell$. 
Left multiplying this equality by $2i\mT\mB(\mA-i\mB^*\mT\mB)^{-1}$
yields an expression for $2i\mT\mB(\bu)$ that can be used in
$\bp+i\mT\mB(\bu) = \bq + 2i\mT\mB(\bu)$ in the last line of
\eqref{DomainDecomposedPb2}. This eventually leads to \eqref{SkeletonEquation}.

\quad\\
Reciprocally if $\bq$ solves \eqref{SkeletonEquation}
and $\bu = (\mA-i\mB^*\mT\mB)^{-1}(\mB^*\bq + \bell)$ and 
$\bp = \bq + i\mT\mB(\bu)$, then we have
$\mA\bu = \mB^*(\bq+i\mT\mB\bu) + \bell = \mB^{*}\bp+\bell$.
On the other hand, using the expression of $\bff$ and
$\mS$, the skeleton equation in \eqref{SkeletonEquation}
writes
\begin{equation*}
  \begin{aligned}
    & \bq + \Pi(\bq + 2i\mT\mB(\mA-i\mB^*\mT\mB)^{-1}(\mB^{*}\bq+\bell)) = 0\\
    & \iff\;\;\bq + \Pi(\bq + 2i\mT\mB(\bu)) = 0\\
    & \iff\;\;\bp-i\mT\mB(\bu) + \Pi(\bp+i\mT\mB(\bu)) = 0
  \end{aligned}
\end{equation*}
This finally proves that the pair $(\bu,\bp)$ satisfies \eqref{DomainDecomposedPb2}
\hfill $\Box$

\quad\\
Next we investigate whether or not the skeleton formulation \eqref{EquivalenceFormulation3}
is uniquely solvable. We will show that this is directly correlated to the unique solvability
of \eqref{PrimaryFormulation}.

\begin{prop}\label{DimensionKernelSkeletonOperator}\quad\\
  Assume \eqref{UnifiedVariationalSetting}-\eqref{SuperBilinearForm}-\eqref{AssumptionAbsorption}-\eqref{LocalUniqueSolvability}.
  The application $(\bv,\bp)\mapsto \bp-i\mT(\bv)$ induces a continuous
  isomorphism from $\mathscr{C}(\mA)\cap \mathscr{X}(\Sigma)$ onto $\mrm{ker}(\Id+\Pi\mS)$.
  As a consequence 
  \begin{equation*}
    \mrm{dim}(\,\mrm{ker}(\Id+\Pi\mS)\,) = \mrm{dim}(\,\mrm{ker}(\mA_{\Omega\times\Gamma})\,).
  \end{equation*}
\end{prop}
\noindent \textbf{Proof:}

First of all, if $(\bv,\bp)\in \mathscr{C}(\mA)\cap \mathscr{X}(\Sigma)$, then
$\bp+i\mT\bv = \mS(\bp-i\mT\bv)$ according to Lemma \ref{CaracScatteringOperator},
and $\bp-i\mT\bv = -\Pi(\bp+i\mT\bv)$ according to \eqref{CaracSingleTraceSpace}.
Combining these two identities leads to $\bp-i\mT\bv\in \mrm{ker}(\Id+\Pi\mS)$.
Next if $(\bv,\bp)\in\mathscr{C}(\mA)\cap \mathscr{X}(\Sigma)$ and
$\bp-i\mT\bv = 0$, then $(\bv,\bp) = (0,0)$ according to Lemma
\ref{EnergyConservationIdentity} hence the injectivity.

Finally if $\bq\in \mrm{ker}(\Id+\Pi\mS)$, then there exists
$(\bv,\bp)\in\mathscr{C}(\mA)$ unique such that $\bp-i\mT\bv = \bq$
according to  Lemma \ref{EnergyConservationIdentity}, and
applying \eqref{DefScatteringOperator}, we obtain $\mS(\bq) =
\mS(\bp-i\mT\bv) = \bp+i\mT\bv$. From this later identity
and $(\Id+\Pi\mS)\bq = 0$ leads to $-\bp+i\mT\bv = \Pi(\bp+i\mT\bv)$
which implies $(\bv,\bp)\in\mathscr{X}(\Sigma)$ according to
Proposition \ref{NonlocalTrCond}. Hence we conclude
$(\bv,\bp)\in\mathscr{C}(\mA)\cap \mathscr{X}(\Sigma)$. \hfill $\Box$

\begin{prop}\label{ClosedRangeSkeleton}\quad\\
  Assume \eqref{UnifiedVariationalSetting}-\eqref{SuperBilinearForm}-\eqref{AssumptionAbsorption}-\eqref{LocalUniqueSolvability}.
  The subspace $\mrm{range}(\Id+\Pi\mS)$  is closed in $\mbH(\Sigma)^*$.
\end{prop}
\noindent \textbf{Proof:}

Define $\Theta: \mbH(\Sigma)^{*}\to \mathscr{H}(\Sigma)$ by
$\Theta(\bq) := (i\mT^{-1}(\bq), \bq)$, which satisfies
$2\Vert \bq\Vert_{\mT^{-1}}^{2} = \Vert \Theta(\bq)\Vert_{\mT\times\mT^{-1}}^2$
for all $\bq\in \mbH(\Sigma)^*$. Taking account that
$\mathscr{C}(\mA)+\mathscr{X}(\Sigma)$ is closed, see
Proposition \ref{ClosedImage}, we are going to prove that
\begin{equation*}
  \mrm{range}(\Id + \Pi\mS) = \Theta^{-1}(\mathscr{C}(\mA)+\mathscr{X}(\Sigma)).
\end{equation*}
Take any $\bp\in \mrm{range}(\Id + \Pi\mS)$. Applying Lemma \ref{EnergyConservationIdentity},
there exists a unique $(\bv,\bq)\in\mathscr{C}(\mA)$ such that $2\bp = (\Id + \Pi\mS)(\bq-i\mT\bv)$.
Since $\mS(\bq-i\mT\bv) = \bq+i\mT\bv$ according to Proposition \ref{DefScatOp2}, and writing
$2\bp = (\Id + \Pi)\bp + (\Id-\Pi)\bp$, we obtain
\begin{equation*}
  \begin{aligned}
    &&(\Id+\Pi)\bp +  (\Id-\Pi)\bp
    & =  \bq-i\mT\bv + \Pi(\bq+i\mT\bv)\\
    \iff&&(\Id+\Pi)\bp +  (\Id-\Pi)\bp
    & =  (\Id+\Pi)\bq - (\Id -\Pi)(i\mT\bv)\\
    \iff&& (\Id+\Pi)(\bp-\bq)
    & =  - (\Id -\Pi)(\bp+i\mT\bv).\\
  \end{aligned}
\end{equation*}
As $(\Id\pm\Pi)/2$ are two mutually orthogonal projectors, see Proposition
\ref{NonlocalTrCond}, we deduce on the one hand that $(\Id+\Pi)(\bp-\bq) = 0$
and $(\Id-\Pi)(\bp+i\mT\bv) = 0$. This eventually leads to 
$\bp-\bq\in\mbX(\Sigma)^\circ$ and $\bp+i\mT\bv\in\mT(\mbX(\Sigma))
\iff i\mT^{-1}\bp - \bv\in \mbX(\Sigma)$. We conclude that
$\Theta(\bp) - (\bv,\bq)\in \mathscr{X}(\Sigma)$.
Hence $\Theta(\bp)\in \mathscr{C}(\mA)+\mathscr{X}(\Sigma)$.

\quad\\
Reciprocally pick an arbitrary $\bp\in \Theta^{-1}(\mathscr{C}(\mA)+\mathscr{X}(\Sigma))$.
This means that $\Theta(\bp) - (\bv,\bq)\in \mathscr{X}(\Sigma)$ for some
$(\bv,\bq)\in\mathscr{C}(\mA)$. As a consequence $(\Id-\Pi)(\bp+i\mT\bv) = 0$
and $(\Id+\Pi)(\bp-\bq) = 0$. Adding these two equations, and taking account that
$\bq+i\mT\bv = \mS(\bq-i\mT\bv)$ according to \eqref{DefScatteringOperator}, 
leads to 
\begin{equation*}
  \begin{aligned}
    &&(\Id+\Pi)(\bp-\bq) &=- 
    (\Id-\Pi)(\bp+i\mT\bv)\\
    \iff
    &&(\Id+\Pi)\bp +  (\Id-\Pi)\bp
    & =  \bq-i\mT\bv + \Pi(\bq+i\mT\bv)\\
    \iff&&
    \bp & = (\Id + \Pi\mS)(\bq-i\mT\bv).
  \end{aligned}
\end{equation*}
\hfill $\Box$

\begin{prop}\label{SameCodimension}\quad\\
  Assume \eqref{UnifiedVariationalSetting}-\eqref{SuperBilinearForm}-\eqref{AssumptionAbsorption}-\eqref{LocalUniqueSolvability}.
  Then
  \begin{equation*}
    \mrm{codim}(\,\mrm{range}(\Id+\Pi\mS)\,) = \mrm{codim}(\,\mrm{range}(\mA_{\Omega\times\Gamma})\,).
  \end{equation*}
\end{prop}
\noindent \textbf{Proof:}

Since $\mrm{range}(\Id+\Pi\mS)$ is closed according to
Proposition \ref{ClosedRangeSkeleton}, we deduce that
$\mrm{codim}(\,\mrm{range}(\Id+\Pi\mS)\,) = \mrm{dim}(\,
\mrm{ker}( (\Id + \Pi\mS)^*)\,)$. Proposition \ref{NonlocalTrCond}, 
in particular the characterization of $\mQ = (\Id +\Pi)/2$ as a $\mT^{-1}$-orthogonal projection,
show that $\Pi^{2} = \Id$ and $\Pi^* = \mT^{-1}\Pi\mT$, so we have
\begin{equation*}
  (\Id + \Pi\mS)^* =
  (\mT\Pi^{*})^{-1}(\Id +\Pi\mT\mS^*\mT^{-1})\mT\Pi^{*}.
\end{equation*}
Setting $\tilde{\mS} := \mT\mS^*\mT^{-1}$, and noting
that $\mT\Pi^{*}:\mbH(\Sigma)\to \mbH(\Sigma)^{*}$ is
an isomorphism, we have
$\mrm{dim}(\,\mrm{ker}( (\Id + \Pi\mS)^*)\,) =
\mrm{dim}(\,\mrm{ker}( \Id + \Pi\tilde{\mS})\,)$.
Let us have a close look at $\tilde{\mS}$, taking
account of the formulas given by Proposition \ref{DefScatOp2}.
Since $\mT^* = \mT$, we obtain 
\begin{equation*}
  \tilde{\mS}=\Id +2i\mT\mB(\mA^*-i\mB^*\mT\mB)^{-1}\mB^*.
\end{equation*}
We see that $\tilde{\mS}$ differs from $\mS$ only in that $\mA$
is replaced by $\mA^{*}$. As a consequence, we can apply
Proposition \ref{DimensionKernelSkeletonOperator}, replacing
$\mA_{\Omega\times\Gamma}$ with $\mA_{\Omega\times\Gamma}^{*}$.
Using \eqref{FredholmAlternative}, this yields
$\mrm{dim}(\,\mrm{ker}( \Id + \Pi\tilde{\mS})\,) =
\mrm{dim}(\,\mrm{ker}(\mA_{\Omega\times\Gamma}^{*})\,) =
\mrm{codim}(\,\mrm{range}(\mA_{\Omega\times\Gamma})\,)$.
\hfill $\Box$

\quad\\
If $\mV_1,\mV_2$ are Banach spaces, a bounded linear map  $\mL:\mV_1\to \mV_2$ is
of Fredholm type if and only if $\mrm{range}(\mL)$ is closed in $\mV_2$,
$\mrm{dim}(\,\mrm{ker}(\mL)\,)<\infty$ and $\mrm{codim}(\,\mrm{range}(\mL)\,)<\infty$.
In this case the index of $\mL$ is the number
$\mrm{index}(\mL):= \mrm{dim}(\,\mrm{ker}(\mL)\,) - \mrm{codim}(\,\mrm{range}(\mL)\,)$.
The results of the present paragraph (in particular Proposition
\ref{DimensionKernelSkeletonOperator}, \ref{ClosedRangeSkeleton} and
\ref{SameCodimension}) lead to the following corollary.

\begin{cor}\quad\\
  Assume \eqref{UnifiedVariationalSetting}-\eqref{SuperBilinearForm}-\eqref{AssumptionAbsorption}-\eqref{LocalUniqueSolvability}.
  The operator $\mA_{\Omega\times\Gamma}:
  \mH(\Omega\times\Gamma)\to \mH(\Omega\times\Gamma)^{*}$ is of Fredholm
  type if and only if $\Id + \Pi\mS:\mbH(\Sigma)^*\to \mbH(\Sigma)^*$ is
  of Fredholm type and, in this case, both operators have the same index.
\end{cor}

\section{Coercivity estimate}

Now we study quantitatively how the inf-sup constant of $\Id+\Pi\mS$
relates to the inf-sup constant of the operator $\mA_{\Omega\times\Gamma}$. Taking
the cue from \cite[\S8]{claeys2021nonself}, we first establish an intermediate
result. Recall that inf-sup constants are defined according to \eqref{DefInfSupCst}.

\begin{prop}\label{EstimationProjectionOblique}\quad\\
  Assume \eqref{UnifiedVariationalSetting}-\eqref{SuperBilinearForm}-\eqref{AssumptionAbsorption}-\eqref{LocalUniqueSolvability}. Then
  \begin{equation*}
    \begin{aligned}
      &\infsup_{\mH(\Omega\times\Gamma)\to \mH(\Omega\times\Gamma)^{*}}(
      \mA_{\Omega\times\Gamma})
      \leq (1+\Vert \mA\Vert) \inf_{\substack{
          \bu\in \mathscr{C}(\mA)\setminus\{0\}\\
          \bv\in \mathscr{X}(\Sigma)\setminus\{0\}}}
      \frac{\Vert \bu + \bv\Vert_{\mT\times\mT^{-1}}}{
        \Vert\bu\Vert_{\mT\times\mT^{-1}}}\\
      &\text{where}\quad
      \Vert \mA\Vert:=\sup_{\bu,\bv\in \mbH(\Omega)\setminus\{0\}}\frac{\vert\langle \bu,\mA(\bv)\rangle\vert}{
        \Vert \bu\Vert_{\mbH(\Omega)} \Vert \bv\Vert_{\mbH(\Omega)} }.
      \end{aligned}
  \end{equation*}
\end{prop}
\noindent \textbf{Proof:}

In the case where $\mathscr{C}(\mA)\cap \mathscr{X}(\Sigma) \neq \{0\}$,
the inf-sup constant vanishes since $\mrm{ker}(\mA_{\Omega\times\Gamma})\neq \{0\}$ according
to Proposition \ref{CorrespondanceKernelIntersection}. So the estimate is automatically
satisfied in this case. We shall assume $\mathscr{C}(\mA)\cap \mathscr{X}(\Sigma) = \{0\}$.
According to Proposition \ref{CorrespondanceKernelIntersection} this
leads to
\begin{equation}\label{TrivialIntersection}
  \begin{aligned}
    & \mrm{ker}(\mA_{\Omega\times\Gamma})\neq \{0\}\\
    & \alpha:=\infsup_{\mH(\Omega\times\Gamma)\to \mH(\Omega\times\Gamma)^{*}}(
      \mA_{\Omega\times\Gamma})>0.
  \end{aligned}
\end{equation}
Now pick any $\bu\in \mathscr{C}(\mA)\setminus\{0\}$
and any $\bv\in \mathscr{X}(\Sigma)\setminus\{0\}$, and set 
$(\bp_{\dir},\bp_{\neu}):= \bu+\bv \in
\mathscr{H}(\Sigma) = \mbH(\Sigma)\times \mbH(\Sigma)^*$.
The invertibility of $\mA_{\Omega\times\Gamma}$ provides the
existence of a unique $\varphi\in  \mbX(\Omega)$ satisfying
$\langle \mA(\varphi),\bw\rangle =
-\langle \mA\mB^{\dagger}(\bp_{\dir}),\bw\rangle 
+\langle  \bp_{\neu},\mB(\bw)\rangle$ for all $\bw\in \mbX(\Omega)$.
In particular
\begin{equation}\label{Estimation1}
  \alpha\,\Vert \varphi\Vert_{\mbH(\Omega)}\leq
  \Vert \mA\Vert\,\Vert \bp_{\dir}\Vert_{\mT} +
  \Vert \bp_{\neu}\Vert_{\mT^{-1}}.
\end{equation}
Set $\phi  = \varphi+\mB^{\dagger}(\bp_{\dir})$ and $\bu_{\dir} =
\mB(\phi) = \mB(\varphi)+\bp_{\dir}$. By construction, for any
$\bw\in \mbH(\Omega)$ satisfying $\mB(\bw) = 0$ we have
$\langle \mA(\phi),\bw\rangle = \langle\bp_{\neu},\mB(\bw)\rangle  = 0$,
which rewrites $\mA(\phi)\in \mrm{ker}(\mB)^{\circ}$. Applying \textit{i)}
of Lemma \ref{ExtendedPolarity} we have $\mA\phi = \mB^*\bu_{\neu}$ for
some $\bu_{\neu}\in \mbH(\Sigma)^{*}$. This implies in particular
$\bu_{\neu} = (\mB\mB^{\dagger})^*\bu_{\neu} = (\mB^{\dagger})^*\mB^*\bu_{\neu} =
(\mB^{\dagger})^*\mA\phi$. From the previous definitions, and the fact that
$\Vert \mB(\bw)\Vert_{\mT}\leq \Vert \bw\Vert_{\mbH(\Omega)}$ and
$\Vert \mB^{\dagger}(\bq)\Vert_{\mbH(\Omega)} =\Vert \bq\Vert_{\mT}$,
we obtain the estimates
\begin{equation}\label{Estimation2}
  \begin{aligned}
    \Vert \phi\Vert_{\mbH(\Omega)}
    & \leq \Vert \varphi\Vert_{\mbH(\Omega)} + \Vert \bp_{\dir}\Vert_{\mT}\\
    \Vert \bu_{\dir}\Vert_{\mT}
    & \leq \Vert \phi\Vert_{\mbH(\Omega)}\\
    \Vert \bu_{\neu}\Vert_{\mT^{-1}}
    & \leq \Vert \mA\Vert\,\Vert \phi\Vert_{\mbH(\Omega)}.
  \end{aligned}
\end{equation}
We have $\mA\phi = \mB^*\bu_{\neu}$ and $\mB\phi = \bu_{\dir}$ hence
$(\bu_{\dir},\bu_{\neu})\in \mathscr{C}(\mA)$ by construction.
On the other hand we have $\bp_{\dir}-\bu_{\dir} = \mB\varphi\in \mbX(\Sigma)$
since $\varphi\in\mbX(\Omega)$ and, for any $\bw\in \mbX(\Sigma)$ we have
$\mB^{\dagger}(\bw)\in\mbX(\Omega)$ hence 
$\langle \bp_{\neu}-\bu_{\neu},\bw\rangle  =
\langle \mA\phi,\mB^{\dagger}\bw\rangle - \langle \mA\phi,\mB^{\dagger}\bw\rangle = 0$,
which implies that $\bp_{\neu}-\bu_{\neu}\in\mbX(\Sigma)^{\circ}$.
Finally we have shown that $(\bu_{\dir},\bu_{\neu})\in\mathscr{C}(\mA)$
and $(\bp_{\dir},\bp_{\neu}) - (\bu_{\dir},\bu_{\neu})\in \mathscr{X}(\Sigma)$
and, since $\bp = \bu + \bv\in\mathscr{C}(\mA)\oplus
\mathscr{X}(\Sigma)$, we conclude that $\bu = (\bu_{\dir},\bu_{\neu})$. 
There only remains to combine \eqref{Estimation1} and
\eqref{Estimation2} to obtain the desired estimate. \hfill $\Box$

\begin{thm}\label{FinalEstimate}\quad\\
  Assume \eqref{UnifiedVariationalSetting}-\eqref{SuperBilinearForm}-\eqref{AssumptionAbsorption}-\eqref{LocalUniqueSolvability}. Then
  \begin{equation*}
    \infsup_{\mH(\Omega\times\Gamma)\to
      \mH(\Omega\times\Gamma)^{*}}(\mA_{\Omega\times\Gamma})
    \leq (1+\Vert \mA\Vert)
    \infsup_{\mbH(\Sigma)^{*}\to\mbH(\Sigma)^{*}}(\Id + \Pi\mS).
  \end{equation*}
\end{thm}
\noindent \textbf{Proof:}

In the case where $\mrm{ker}(\mA_{\Omega\times\Gamma})\neq \{0\}$
we also have $\mrm{ker}(\Id + \Pi\mS)\neq \{0\}$ according to
Proposition \ref{DimensionKernelSkeletonOperator} and, in this situation, the
desired estimate is satisfied, with both sides of the estimate equal
to $0$. Hence we can assume that 
$\mrm{ker}(\mA_{\Omega\times\Gamma}) = \{0\}$ and in this situation
both $\mA_{\Omega\times\Gamma}:\mH(\Omega\times\Gamma)\to \mH(\Omega\times\Gamma)^*$
and $\Id+\Pi\mS:\mbH(\Sigma)\to\mbH(\Sigma)^*$ are are injective with closed
range. Pick an arbitrary $\bff\in \mbH(\Sigma)^*$. According to Lemma
\ref{EnergyConservationIdentity}, there exists a unique pair
$\bu = (\bu_\dir,\bu_\neu)\in\mathscr{C}(\mA)$ such that $\bff =
\bu_{\neu} - i\mT(\bu_{\dir})$ and we have $\Vert \bff\Vert_{\mT^{-1}}\leq
\sqrt{2}\Vert \bu\Vert_{\mT\times\mT^{-1}}$ which, for $\bff\neq 0$,
re-writes as
\begin{equation*}
  \frac{\Vert \bu\Vert_{\mT\times\mT^{-1}}}{\Vert \bff\Vert_{\mT^{-1}\textcolor{white}{\times\mT}}}\geq \frac{1}{\sqrt{2}}.
\end{equation*}
Next set $\bg = (\Id+\Pi\mS)\bff$ and $\bp = (\bp_{\dir},\bp_{\neu})
= (\mT^{-1}(\bg), -i\bg)/2$. We have in particular $\Vert \bg\Vert_{\mT^{-1}} =
\sqrt{2}\Vert \bp\Vert_{\mT\times\mT^{-1}}$. Since $\mS(\bff) = \mS(\bu_{\neu} - i\mT(\bu_{\dir})) =
\bu_{\neu} + i\mT(\bu_{\dir})$ according to Proposition
\ref{DefScatOp2}, we obtain 
\begin{equation*}
  \begin{aligned}
    \bu_{\neu} - i\mT(\bu_{\dir})+\Pi(\bu_{\neu} + i\mT(\bu_{\dir}))
    & = \bff + \Pi\mS(\bff) \\
    & = \bg = (\Id+\Pi)\bg/2 + (\Id-\Pi)\bg/2\\
    & = (\Id+\Pi)\bp_\neu -i (\Id-\Pi)\mT(\bp_\dir)\\
    & = \bp_\neu-i\mT(\bp_{\dir}) + \Pi(\bp_\neu+i\mT(\bp_\dir))
  \end{aligned}
\end{equation*}
Re-arranging the terms in the equality above so as to move all contributions
involving $\Pi$ in the right hand side, we obtain 
$-(\bp_\neu-\bu_{\neu})+i\mT(\bp_{\dir}-\bu_{\dir}) =
\Pi( (\bp_\neu-\bu_{\neu}) + i\mT(\bp_{\dir}-\bu_{\dir}))$.
According to Proposition \ref{NonlocalTrCond}, this implies that
$(\bp_\dir,\bp_\neu) - (\bu_\dir,\bu_\neu)\in \mathscr{X}(\Sigma)$.
Since we have $(\bu_\dir,\bu_\neu)\in\mathscr{C}(\mA)$ by construction,
we can apply Proposition \ref{EstimationProjectionOblique}
which, for $\bff\neq 0$, yields
\begin{equation*}
  \frac{\Vert (\Id+\Pi\mS)\bff\Vert_{\mT^{-1}}}{\Vert\bff\Vert_{\mT^{-1}}} =
  \frac{\Vert\bg\Vert_{\mT^{-1}}}{\Vert\bff\Vert_{\mT^{-1}}}\geq
  \frac{\Vert\bp\Vert_{\mT\times\mT^{-1}}}{\Vert\bu\Vert_{\mT\times\mT^{-1}}}
  \geq \infsup_{\mH(\Omega\times\Gamma)\to
      \mH(\Omega\times\Gamma)^{*}}(\mA_{\Omega\times\Gamma})
    /(1+\Vert \mA\Vert).
\end{equation*}
This establishes the desired estimate, since this holds for any
$\bff\in \mbH(\Sigma)^*\setminus\{0\}$.
\hfill $\Box$

\quad\\
The estimate provided by Theorem \ref{FinalEstimate} is remarkable in several
respects. First of all it holds even if $\mrm{ker}(\mA_{\Omega\times\Gamma})$ is non-trivial.
Secondly it does not involve any hidden ``$C>0$'' constant. In particular it does not
involve any frequency dependency, although the infsup constant of
$\mA_{\Omega\times\Gamma}$ a priori depends itself on the frequency.
This means that, to estimate the frequency dependency of the infsup constant of
$\Id+\Pi\mS$, it suffices to derive such an estimate for $\mA_{\Omega\times\Gamma}$.
A further striking feature is that the number of subdomains $\mJ$
does not come into play in this estimate. 

\quad\\
As an interesting additional result in the perspective of an effective linear solve, 
the contractivity of $\Pi$ and $\mS$ leads to the coercivity of the operator $\Id+\Pi\mS$.
The next result can be combined with Theorem \ref{FinalEstimate} to obtain
an effective estimate of the coercivity constant.

\begin{cor}\label{FinalCoercivityEstimate}\quad\\
  Assume\eqref{UnifiedVariationalSetting}-\eqref{SuperBilinearForm}-\eqref{AssumptionAbsorption}-\eqref{LocalUniqueSolvability}. Then
  $\Id+\Pi\mS:\mbH(\Sigma)^*\to \mbH(\Sigma)^*$ is coercive with
  respect to the scalar product induced by $\mT^{-1}$ and we have
  \begin{equation*}
    \inf_{\bq\in\mbH(\Sigma)^{*}\setminus\{0\}}
    \frac{\Re e\{\langle (\Id + \Pi\mS)\bq,\mT^{-1}\overline{\bq}\rangle\}}{
      \Vert \bq\Vert_{\mT^{-1}}^{2}}\geq \frac{1}{2}\big(
    \infsup_{\mbH(\Sigma)^{*}\to\mbH(\Sigma)^{*}}(\Id + \Pi\mS)\;\big)^{2}.
  \end{equation*}  
\end{cor}
\noindent \textbf{Proof:}

For any $\bq\in \mbH(\Sigma)^{*}\setminus\{0\}$,
\begin{equation*}
  \begin{aligned}
    & \Vert \bq\Vert_{\mT^{-1}}^{2} \geq \Vert \Pi\mS(\bq)\Vert_{\mT^{-1}}^{2}
    = \Vert (\Id+\Pi\mS)\bq -  \bq\Vert_{\mT^{-1}}^{2}\\
    &
    \textcolor{white}{\Vert \bq\Vert_{\mT^{-1}}^{2} \geq \Vert \Pi\mS(\bq)\Vert_{\mT^{-1}}^{2}}
    = \Vert (\Id+\Pi\mS)\bq\Vert_{\mT^{-1}}^{2} + \Vert\bq\Vert_{\mT^{-1}}^{2}
    - 2\Re e \{ \langle (\Id + \Pi\mS)\bq,\mT^{-1}\overline{\bq}\rangle\}\\
    & \Longrightarrow\quad
    \Re e \{ \langle (\Id + \Pi\mS)\bq,\mT^{-1}\overline{\bq}\rangle\}/\Vert \bq\Vert_{\mT^{-1}}^{2}
    \geq \big(\Vert (\Id+\Pi\mS)\bq\Vert_{\mT^{-1}}/\Vert \bq\Vert_{\mT^{-1}}\big)^{2}/2.
  \end{aligned}
\end{equation*}
\hfill $\Box$

\quad\\
We conclude this article illustrating how the previous results lead to
estimations of the coercivity constant of the skeleton operator for a
concrete case.

\begin{example}\quad\\
  Consider the case $\RR^d = \RR^2$ or $\RR^3$.  Assume that $\mu = 1$, $\kappa = k\in (0,+\infty)$, and
  choose $\mA_{\Gamma}$ as in Example \ref{RobinBC1} with $\langle \Lambda(u), v\rangle = k\int_{\Gamma} u v d\sigma$
  which models the Robin condition $\partial_{\bn}u -ik u = 0$ on $\Gamma$. So we.
  Assume in addition that $\Omega$ is a convex polyhedron. Then we have
  \begin{equation*}  
    \langle \mA_{\Omega\times\Gamma}(u,p),(v,q)\rangle = \int_{\Omega}\nabla u\nabla v - k^2 uv d\bx - ik\int_{\Gamma}uv d\sigma
    +\int_{\Gamma}q\mT_{\Gamma} p \,d\sigma.
  \end{equation*}
  Let us take $\gamma = 1/k$ for the parameter involved in \eqref{H1fctspace}.
  From these choices, and proceeding like in \cite[Lem.2.4]{zbMATH07248609} for dealing
  with boundary terms on $\Gamma$,  we see that the continuity modulus $\Vert\mA\Vert$ (as defined in
  Proposition \ref{EstimationProjectionOblique}) can be bounded independently of $k$. On the other hand,
  we know from \cite{MR2692949} that
  \begin{equation*}
    \infsup_{\mH(\Omega\times\Gamma)\to \mH(\Omega\times\Gamma)^*}(\mA_{\Omega\times\Gamma})
    \geq \mathop{\mathcal{O}}_{k\to \infty}(1/k).
  \end{equation*}
  We can now plug this estimate into Theorem \ref{FinalEstimate}, and we see that the inf-sup constant
  of $\Id+\Pi\mS$ admits also a lower bound that behaves like $\mathcal{O}(1/k)$ for $k\to \infty$.
  Finally combining with Corollary \ref{FinalCoercivityEstimate}, we see that the coercivity constant of the skeleton formulation
  behaves like $\mathcal{O}(1/k^2)$ i.e.
  \begin{equation*}
    \inf_{\bq\in\mbH(\Sigma)^{*}\setminus\{0\}}
    \Re e\{\langle (\Id + \Pi\mS)\bq,\mT^{-1}\overline{\bq}\rangle\}/
      \Vert \bq\Vert_{\mT^{-1}}^{2}\geq  \mathop{\mathcal{O}}_{k\to \infty}(1/k^2).
  \end{equation*}

\end{example}


\begin{thebibliography}{10}

\bibitem{Bendali2006}
A.~Bendali and Y.~Boubendir.
\newblock Non-overlapping domain decomposition method for a nodal finite
  element method.
\newblock {\em Numerische Mathematik}, 103(4):515--537, Jun 2006.

\bibitem{MR2759829}
H.~Brezis.
\newblock {\em Functional analysis, {S}obolev spaces and partial differential
  equations}.
\newblock Universitext. Springer, New York, 2011.

\bibitem{MR1618464}
O.~Cessenat and B.~Despres.
\newblock Application of an ultra weak variational formulation of elliptic
  {PDE}s to the two-dimensional {H}elmholtz problem.
\newblock {\em SIAM J. Numer. Anal.}, 35(1):255--299, 1998.

\bibitem{zbMATH00195024}
P.G. Ciarlet.
\newblock {\em Introduction to numerical linear algebra and optimization.}
\newblock Camb. Texts Appl. Math. Cambridge etc.: Cambridge University Press,
  1988.

\bibitem{claeys2019new}
X.~Claeys.
\newblock Non-local variant of the {O}ptimised {S}chwarz {M}ethod for arbitrary
  non-overlapping subdomain partitions.
\newblock {\em ESAIM: M2AN}, 55(2):429--448, 2021.

\bibitem{claeys2021nonself}
X.~Claeys.
\newblock {Nonselfadjoint impedance in Generalized Optimized Schwarz Methods}.
\newblock {\em IMA Journal of Numerical Analysis}, November 2022.

\bibitem{MR4507159}
X.~Claeys, F.~Collino, and E.~Parolin.
\newblock Nonlocal optimized schwarz methods for time-harmonic
  electromagnetics.
\newblock {\em Adv. Comput. Math.}, 48(6):Paper No. 72, 2022.

\bibitem{MR4433119}
X.~Claeys and E.~Parolin.
\newblock Robust treatment of cross-points in optimized {S}chwarz methods.
\newblock {\em Numer. Math.}, 151(2):405--442, 2022.

\bibitem{MR1764190}
F.~Collino, S.~Ghanemi, and P.~Joly.
\newblock Domain decomposition method for harmonic wave propagation: a general
  presentation.
\newblock {\em Computer Methods in Applied Mechanics and Engineering},
  184(2):171 -- 211, 2000.

\bibitem{MR1291197}
B.~Despr\'{e}s.
\newblock {\em M\'{e}thodes de d\'{e}composition de domaine pour les
  probl\`emes de propagation d'ondes en r\'{e}gime harmonique. {L}e
  th\'{e}or\`eme de {B}org pour l'\'{e}quation de {H}ill vectorielle}.
\newblock Institut National de Recherche en Informatique et en Automatique
  (INRIA), Rocquencourt, 1991.
\newblock Th\`ese, Universit\'{e} de Paris IX (Dauphine), Paris, 1991.

\bibitem{MR4480644}
B.~Despr\'{e}s, A.~Nicolopoulos, and B.~Thierry.
\newblock Optimized transmission conditions in domain decomposition methods
  with cross-points for {H}elmholtz equation.
\newblock {\em SIAM J. Numer. Anal.}, 60(5):2482--2507, 2022.

\bibitem{Gander2013}
M.~Gander and F.~Kwok.
\newblock On the applicability of {L}ions' energy estimates in the analysis of
  discrete optimized schwarz methods with cross points.
\newblock {\em Lecture Notes in Computational Science and Engineering}, 91, 01
  2013.

\bibitem{MR3519297}
M.J. Gander and K.~Santugini.
\newblock Cross-points in domain decomposition methods with a finite element
  discretization.
\newblock {\em Electron. Trans. Numer. Anal.}, 45:219--240, 2016.

\bibitem{zbMATH07020343}
M.J. {Gander} and H.~{Zhang}.
\newblock {A class of iterative solvers for the Helmholtz equation:
  factorizations, sweeping preconditioners, source transfer, single layer
  potentials, polarized traces, and optimized Schwarz methods}.
\newblock {\em {SIAM Rev.}}, 61(1):3--76, 2019.

\bibitem{zbMATH07248609}
I.G. {Graham}, E.A. {Spence}, and J.{Zou}.
\newblock {Domain decomposition with local impedance conditions for the
  Helmholtz equation with absorption}.
\newblock {\em {SIAM J. Numer. Anal.}}, 58(5):2515--2543, 2020.

\bibitem{MR1335452}
T.~Kato.
\newblock {\em Perturbation theory for linear operators}.
\newblock Classics in Mathematics. Springer-Verlag, Berlin, 1995.
\newblock Reprint of the 1980 edition.

\bibitem{zbMATH01446717}
W.~McLean.
\newblock {\em Strongly elliptic systems and boundary integral equations}.
\newblock Cambridge: Cambridge University Press, 2000.

\bibitem{MR2692949}
J.~M. Melenk.
\newblock {\em On generalized finite-element methods}.
\newblock ProQuest LLC, Ann Arbor, MI, 1995.
\newblock Thesis (Ph.D.)--University of Maryland, College Park.

\bibitem{Modave2020b}
A.~{Modave}, A.~{Royer}, X.~{Antoine}, and C.~{Geuzaine}.
\newblock {A non-overlapping domain decomposition method with high-order
  transmission conditions and cross-point treatment for Helmholtz problems}.
\newblock {\em {Comput. Methods Appl. Mech. Eng.}}, 368:23, 2020.
\newblock Id/No 113162.

\bibitem{parolin:tel-03118712}
E.~Parolin.
\newblock {\em {Non-overlapping domain decomposition methods with non-local
  transmissionoperators for harmonic wave propagation problems}}.
\newblock Theses, {Institut Polytechnique de Paris}, December 2020.

\bibitem{MR3013465}
C.~Pechstein.
\newblock {\em Finite and boundary element tearing and interconnecting solvers
  for multiscale problems}, volume~90 of {\em Lecture Notes in Computational
  Science and Engineering}.
\newblock Springer, Heidelberg, 2013.

\bibitem{zbMATH01022519}
W.~Rudin.
\newblock {\em {Functional analysis. 2nd ed}}.
\newblock New York, NY: McGraw-Hill, 2nd ed. edition, 1991.

\bibitem{MR2361676}
O.~Steinbach.
\newblock {\em Numerical approximation methods for elliptic boundary value
  problems}.
\newblock Springer, New York, 2008.
\newblock Finite and boundary elements, Translated from the 2003 German
  original.

\bibitem{zbMATH04136456}
T.~von Petersdorff.
\newblock Boundary integral equations for mixed {Dirichlet}, {Neumann} and
  transmission problems.
\newblock {\em Math. Methods Appl. Sci.}, 11(2):185--213, 1989.

\end{thebibliography}

\end{document}